\input amstex
\documentstyle{amsppt}
\magnification=\magstep1
\pageheight{8.5truein}
\pagewidth{5.5truein}
\hcorrection{.7in}
\parindent=10pt
\nologo

\NoRunningHeads

\define\go{{\frak g}^{*}\times\Omega}

\define\pgo{\text{Prim}(C^*(G,\Omega))}

\define\twoskip{\smallskip\smallskip}
\define\dd{\text{d}}
\magnification=1100
\font\fod=cmr10\magnification=1200

\overfullrule=0pt

\document

\baselineskip=20pt

\raggedbottom

\topmatter
\title
{Kirillov theory for $C^*(G,\Omega)$}
\endtitle

\author
Dean Moore
\endauthor

\abstract

Let $G$ be a simply connected
nilpotent Lie group with Lie algebra $\frak g$;
let $\frak g^*$ be the dual of $\frak g$.
Let $\Omega$ be a locally compact second countable
Hausdorff space with a continuous $G$ action, and
let $C^*(G,\Omega)$ be the corresponding transformation group $C^*$
algebra.
We construct
a continuous surjective map $\phi$ from a quotient space, $\frak g^*\times\Omega/\sim$,
which is a
homeomorphism from $\frak g^*\times\Omega/\sim$ to Prim$(C^*(G,\Omega))$.

We also describe a character theory for $C^*(G,\Omega)$ which
generalizes Kirillov character theory for $G$.

\endabstract
\endtopmatter

\topmatter
\title\chapter{1} Introduction\endtitle
\endtopmatter

\twoskip

Our primary references for this paper are
Corwin and Greenleaf \cite{3}; fundamental ideas were provided
by Siegfried Echterhoff's paper \cite{7},
and the original inspiration was Dana P. Williams's 
paper \cite{26}.

\smallskip

Let $G$ be a locally compact group acting on a second countable
locally compact Hausdorff space $\Omega$.

A basic question of $C^*$ theory is an
explicit description of Prim$(C^*(G,\Omega))$. For $G$ a connected,
simply connected nilpotent Lie group and $\Omega$ a point, a
description of Prim$(C^*(G,\Omega))$ as a set
was given in the classic paper 
of Kirillov \cite{21} where he showed that the natural map from
$\frak g^*/G\mapsto\widehat{G}$ was continuous with $\frak g^*/G$ in
the quotient topology. 
Kirillov conjectured that this was a homeomorphism; this was
first proved by
Brown \cite{1}. When $\Omega$ is not a
single point, but $G$ is abelian, a complete description of
Prim$(C^*(G,\Omega))$, including a natural description of the topology,
was provided by
Dana P. Williams \cite{26}. In \cite{7}, for general $G$, Siegfried Echterhoff gave a complete
description but required strong hypotheses for the action.


 
We define a map 
$\phi:\go/\sim\mapsto\space \pgo$, and, with a mild restriction upon the action,
show that $\phi$ is a continuous open surjective map when $\go$ is given
the product topology. We further provide a homeomorphism $\psi$
from a quotient space $\sim$ of $\go$ to $\pgo$.

When $C^*(G,\Omega)$ is Type I, we provide a character theory
corresponding to the Kirillov character theory of $G$.
\twoskip\twoskip\twoskip
\topmatter
\title Sections\endtitle
\endtopmatter

This paper is divided into three sections: Section one consists of
preliminary
results as well as setting notation. In section two we prove the
main results of this paper,

 \ \ \ \ 1. An explicit parameterization of $\pgo$, {\it and}, 

\ \ \ \ 2. a computation of the topology of $\pgo$. 

\parindent=10pt
 In section three, we generalize the Kirillov
character formula for $C^*(G,\Omega)$ under stronger hypotheses of the
action. 

\twoskip\twoskip\twoskip
\topmatter
\title Credits\endtitle
\endtopmatter

This paper constitutes part of the author's doctoral thesis at the
the University of Colorado at Boulder, as well as abstraction done 
later. I would like to thank my advisor, Jeff
Fox, for his help.
I would also like to extend a word of thanks to Larry Baggett, Carla
Farsi, Arlan Ramsay, and Marty Walter for helpful discussions. 
Another
word of thanks is extended to Judy Packer, my first doctoral
advisor prior to Jeff Fox, for her suggestions and words of 
encouragement.

\vfill\eject

\noindent{\fod\bf Section 1, Preliminaries}
\twoskip

\twoskip
\noindent{\fod\bf Section 1.1}
\twoskip
\noindent{\fod\bf $C^*$ algebras and representation spaces}

\twoskip
\noindent{\bf Definition 1.}
    
Denote by $\Omega/G$ the orbit space with the quotient
topology.

For $\phi\in C_0(\Omega)$, define ${}^s\phi$ 
by
${}^s\phi(x) = \phi(s^{-1}\cdot x)$.

Denote by $\Cal K(G)$ the space of closed subgroups of $G$, given
the relativized compact-open topology of Fell \cite{10}.

Let $\Cal Q(G)=\{\langle H,T\rangle\ \big |\ H\in\Cal K(G)\ ,\ T\in\text{Rep}(G)\}$.

Let $x\in\Omega$; by $G_x$ we denote the stabilizer of $x$ in $G$.
These are assumed connected.

For $x\in\Omega$, $f\in\frak g^*$,
by $\frak p_x$ we will denote a polarizing subalgebra of
$\frak g_x$ for $f$. Let $P_x=\text{exp}(\frak p_x)$.

Denote by $\frak p$ 
an isotropic subalgebra, not necessarily polarizing.
Let $P=\text{exp}(\frak p)$.

Also:

Define $\chi_{f,P_x}$ to be the character 
of $P_x=\text{exp}(\frak p_x)$ given by $\chi_{f,P_x}(\text{exp}(X))=
e^{i\cdot f(X)}$.

\item{}{Let $\tau_{{}_{f,x}}=\text{ind}_{P_x}^{G_x}(\chi_{{}_{f,P_x}})$ be
an irreducible representation of $G_x$,}
\item{}{Let $\tau_{{}_{f,x}}^\prime=\text{ind}_P^{G_x} be
(\chi_{{}_{f,P}})$, a representation of $G_x$, not in general irreducible.}

We assume representations non-degenerate.

\twoskip

\noindent{\bf Definition 2.}
      A {\it covariant representation} $L=(V_L, M_L)$ of $(G,\Omega)$
on a Hilbert space $H_L$ consists of a uniformly bounded strongly
continuous unitary representation $V$ of $G$ on $H_L$, and a 
norm-decreasing non-degenerate *-preserving representation of
$C_0(\Omega)$, $M$ on $H_L$ such that
$V(s)M(\phi) V(s^{-1})=M({}^s\phi)$.

 Let $A$ be a fixed $C^*$ algebra.
We use the hull-kernel topology on $\widehat A$; see \cite{24}, 
Theorem 5.4.6. 
\smallskip
On the set $\Cal I(A)$ of closed ideals of $A$ we
 use the topology developed by Dana P. Williams on pg. 338 \cite{26}.
This is the topology that has as its sub-base the sets 
$
\{\Cal O_J\}_{{}_{J\in\Cal I(A)}},\ \text{where}\ 
\Cal O_J=\{I\in\Cal I(A)\ |\ I\not\supset J\}.
$
 On Prim$(A)$ this topology restricts to the usual
hull-kernel (Jacobson) topology. One can see that this topology is almost Fell's
\lq\lq inner
hull kernel" topology; see \cite{10}.
\twoskip

\noindent {\bf Definition 3.} Let $\pi$ and $\rho$
be representations of a
locally compact group $G$.

The representation $\rho$ is {\it weakly contained}
in $\pi$ if every positive definite matrix coefficient 
$\langle \rho(\xi),\xi\rangle$ can be approximated uniformly on
compacta of $G$ by finite sums of positive definite matrix coefficients
$\langle \pi(x)\xi^\prime,\xi^\prime\rangle$.
For $\rho$ weakly contained in
$\pi$, we use the notation $\rho\prec\pi$.
The {\it spectrum} of $\pi$ is the set of all representations
weakly contained in $\pi$. 
 
Let $H\subseteq G$ be a closed subgroup, and
let $f_0$ be a non-negative, real-valued function in
     $C_c(G)$ that does not vanish at the identity element. For the
     remainder of this paper, let $\mu_{{}_H}$ be the left Haar 
measure
     on $H$ defined by
$$
\int_Hf_0(t)\dd\mu_{{}_H}(t)=1.
$$
\noindent Such a choice is referred to as a continuous (\lq\lq 
smooth")
choice of Haar
measures, and has the property that $H\mapsto\int_Hf\dd\mu_{{}_H}$ is
continuous on
$\Cal K(G)$ for each $f\in C_c(G)$; see the theorem of \cite{15}, page 908.

\twoskip
\proclaim{Lemma 1} Assume $\{f_n\}_{n=1}^\infty\subseteq 
C_c(G)$ converges to
$f\in C_c(G)$ in the inductive limit topology and $H_n\to H$ in $\Cal
K(G)$. Then

$$
\int_{H_n}f_n\dd\mu_{{}_{H_n}}\longrightarrow\int_{{}_H}f
\dd\mu_{{}_H}.
$$
\endproclaim

\demo{Proof}\enddemo Follows by an $\frac{\epsilon}{2}$ argument.\qed

\twoskip
\noindent{\bf Definition 4.}
      Now define $\Cal E(X)$ to be the union of all the spaces $C(A)$, where
$A\in\Cal K(S)$, thus the elements of $\Cal E(X)$ are complex-valued functions
$f$ such that for the domain of $f$, $D(f)$, we have $D(f)\in\Cal K(S)$. 
Identifying a function with its graph we have $\Cal E(X)\subseteq\Cal K(S)$. We
will always consider the set $\Cal E(X)$ as being equipped with the
(relativized) semicompact-open topology.
\twoskip

\twoskip
\noindent{\bf Definition 5.} Let $X$ be an arbitrary locally compact
space. Let $\Cal K(X)$ denote the
set of all closed subsets of $X$, and equip $\Cal K(X)$ with the topology
whose open sets have as basic open neighborhoods 
$$
U(C,\Cal F)=\{F\in \Cal K(X)\ |\ F\cap C=\emptyset,\  F\cap O\ne\emptyset\ 
\text{for all}\ O\in\Cal F\}$$
where $C$ is a compact subset of $X$ and 
$\Cal F$ is a finite family of open sets of $X$. This topology is call the
{\it compact-open topology} of $\Cal K(X)$. We note $\Cal K(X)$ equipped with
the compact-open topology
is a compact Hausdorff space, see \cite{10}, Theorem 1.\qed

\twoskip
\proclaim{Lemma 2} Let $\{\langle K_i,T_i\rangle\}$ be a net of
elements of $\Cal Q(G)$ and $\langle K,T\rangle$ an element of $\Cal Q(G)$. 
Then $\langle K_i,T_i\rangle\to\langle K,T\rangle$ if and only if, for
each finite sequence
$\phi_1,\cdots,\phi_n$ of functions of positive type on $K$ associated
with $\langle K,T\rangle$, and each subnet of $\{\langle K_i,T_i\rangle\}$,
there exists
\roster
\item a
subnet $\{\langle K^{\prime j},T^{\prime j}\rangle\}$ of that subnet,
and
\item for
each $j$ and each $r=1,\cdots,n$ a finite sum $\phi_r^j$ of functions of
positive type associated with $\langle K^{\prime j},T^{\prime j}\rangle$ 
such that $\phi_r^i\to\phi_r$ in $\Cal E_s(G)$ for each $r$.
\endroster\endproclaim
\demo{Proof}\enddemo See \cite{13} Theorem 3.1', page 439.\qed
\smallskip

\twoskip

\noindent{\bf\fod\bf Section 1.2}

\twoskip
\noindent{\fod\bf Nilpotent Lie algebras and groups, and their representation spaces}
\twoskip

The purpose of this section is to give some basic information
that is needed for this paper
about nilpotent Lie algebras, groups, and the representation space.

\twoskip

     Let $G$ denote a Lie group and $\frak g$ its Lie algebra.

\twoskip
\noindent{\bf Definition 6.}
     The {\it adjoint representation}, ad of $\frak g$ on
$\frak g$ is defined as $\text{ad}_x:\frak g\mapsto \text{GL}(\frak g)$ by
$\text{ad}_x(y)=[x,y]$, for all $y\in\frak g$ (here $[\cdot ,\cdot]$
denotes Lie bracket on $\frak g$).\smallskip

\twoskip
\noindent{\bf Definition 7.}
     The Lie algebra $\frak g$ is said to be {\it nilpotent} if
$\text{ad}_x$ is
a nilpotent endomorphism of $\frak g$, for all $x\in\frak g$.

     The Lie group $G$ is {\it nilpotent} if $\frak g$ is nilpotent.

\twoskip

Now we briefly describe the representation theory
of nilpotent Lie groups that we use. We follow the notation
of Corwin and Greenleaf \cite{3} and use it as our main reference.

\twoskip

     Let $G$ be a connected, simply-connected nilpotent Lie group.
     Denote the dual space of $\frak g$ by $\frak g^*$. 

     $G$ acts on
     $\frak g^*$ by the {\it coadjoint map}, $\text{Ad}^*$:
     for $x\in G$, $Y\in\frak g$, and $f\in\frak g^*$, define
$
(\text{Ad}^*(x)f)(Y)=f(\text{Ad}(x^{-1})Y),
\ Y\in\frak g,\ f\in\frak g^*,\
x\in G.
$

If $f\in\frak g^*$. Define the {\it coadjoint orbit} of $f$
in $\frak g^*$ to be
$\text{Ad}^*(G)f$.

Let $f\in \frak g^*$. A subspace $\frak p\subseteq\frak g$ 
is called {\it isotropic} 
if $f|_{\frak p}= 0$. If $\frak p$ is a maximally isotropic subspace
of $\frak g$ which is also a subalgebra, then $\frak p$ is called a
{\it polarization}, or a 
 {\it maximal subordinate subalgebra} for $f$.

For any $X, Y\in\frak g$, we use the Campbell-Baker-Hausdorff
formula to form $\text{exp}(X)*\text{exp}(Y)$; see Corwin and Greenleaf \cite{3}, pg. 11.

 \twoskip
We have the following result, which seems well known:

\twoskip
\proclaim{Lemma 3} Let $\frak p$ be a polarizing subalgebra of
$\frak g$ for $f\in\frak g^*$.
We may define a one-dimensional representation $\chi_{f,P}$ of
 $P=\text{exp}(\frak p)$ by
$\chi_{f,P}(\text{exp}(X)) = e^{i\cdot f(X)}$.
\endproclaim

\twoskip

We may induce the representation $\chi_{f,P}$ from
$P$ to a representation $\pi_{f,P}$ of $G$; for 
details on induced representations, see \cite{22} and \cite{23}.

\smallskip

Give $\frak g^*/G$ the quotient topology and $\widehat{G}$ the hull-kernel
topology. By \cite{1} and  \cite{21}  we have the following important theorem:

\twoskip
\proclaim{Theorem 1}

\smallskip
\item{(1)} Let $f\in\frak g^*$, and $\frak p$ be a polarizing subalgebra
for $f$.
Let $\pi_{f,P} = \text{ind}_P^G(\chi_{f,P})$.
The representation $\pi_{f,P}$
is irreducible, and up to equivalence,
 every irreducible representation of $G$ is obtained this way.

\smallskip

\item{(1)} We have $\pi_{f,P}\cong\pi_{f^\prime,P^\prime}\ \iff\ \ \exists
\ g\in G$ such that $\text{Ad}^*(g)f = f^\prime$.

\item{(3)} The induced map (the Kirillov map) 
$\frak g^*/Ad^*(G)\mapsto \hat G$ is a homeomorphism.

\endproclaim

\twoskip

Let $G$ be a nilpotent Lie group, $S$ a closed connected subgroup,
having Lie algebras $\frak g$ and $\frak s$, respectively. Define
$\frak s^\perp$ to be the set of linear functionals in $\frak g^*$
that are zero on $\frak s$, i.e.,
$
\frak s^\perp = \{g\in\frak g^*\ |\ g{|_\frak s}=0\}.
$

\twoskip
\proclaim{Lemma 4}
Let $G$ be a simply connected nilpotent Lie group, $S$ a closed 
connected
subgroup;
assume $f\in\frak g^*$ satisfies $f([\frak s,\frak s])=0$ so 
$\chi_{{}_f}(exp(Y))=e^{i\cdot f(Y)}$ is a one dimensional 
representation of
$S$. The representation $W=ind_S^G(\chi_{{}_f})$ weakly contains
$
\bigl\{\pi_{{}_{f^\prime}}\in\widehat G \ |
\ f^\prime\in f+\frak s^\perp\bigr\},
$
\noindent in fact,
$
\text{Sp}(W)=\text{Fell-closure}\bigl(\bigl\{\pi_{f^\prime}\in
\widehat G\ |
\ f^\prime\in f+{\frak s}^\perp\bigr\}\bigr).
$
\endproclaim

\demo{Proof}\enddemo 
A proof may be found in in \cite{4}, Theorem N.2.5.

\twoskip
\noindent {\bf Definition 8.}
Let $G$ a connected, simply connected nilpotent
Lie group. By Kirillov, \cite{21}, the irreducible representations of $G$
are in one-to-one correspondence with $\text{Ad}^*$ orbits in $\frak g^*$.

If $W\in\widehat{G}$, we denote its correspondence orbit by $\Omega_W$.
\twoskip

The following result was proven by Joy \cite{19}. 

\twoskip
\proclaim{Lemma 5 (Joy's Lemma)}
 Let $G$ be a real, connected, simply connected nilpotent
Lie group. Let $\langle H_n,S_n\rangle\to \langle H,S\rangle$ in
$\Cal Q(G)$. If $f\in\frak g^*$ such that $f|_{\frak h}\in\Omega_S$, then
for every sub-sequence of $\{\langle H_n,S_n\rangle\}_{n=1}^\infty$, there is
a sub-sequence $\{\langle H_{n_i},S_{n_i}\rangle\}_{i=1}^\infty$ such that for
each $i$, there exists $f_i\in\frak g^*$ such that 
$f_i|_{\frak h_{n_i}}\in\Omega_{S_{n_i}}$ and $f_i\to f$ in $\frak g^*$.
\endproclaim
\demo{Proof}\enddemo See \cite{19}, page 138.\qed

\twoskip

\noindent{\fod\bf Section 1.3}

\twoskip
\noindent{\fod\bf Induced representations, and ideals of $C^*(G,\Omega)$}

\twoskip

Let $G$ be a group, acting on a locally
compact Hausdorff space $\Omega$.

\twoskip

Let $x\in\Omega$ and $G_x=\{g\in G\ |\ gx = x\}$ be the stability subgroup
of $x$. Let $\rho_x : C_0(\Omega)\mapsto \Bbb C$ be the representation
of $C_0(\Omega)$ given by evaluation at $x$; for $\phi\in C_0(\Omega)$,
$\rho_x(\phi)=\phi(x)$.

Let $\frak g_x$ be the Lie algebra of $G_x$.

Let $\tau$ be a representation of $G_x$ on the Hilbert space $H_{\tau}$.
The pair $(\tau,\rho_x)$ forms a covariant pair for $C^*(G_x,\Omega)$
with the representation of $C^*(G_x,\Omega)$ defined as follows:
for $v\in H_{\tau}$, $\phi\in C_0(\Omega)$, $x\in\Omega$, $g\in G_x$, define
$
\rho_x(\phi)(v)=\phi(x)\cdot v.
$

We may induce the representation $(\tau,\rho_x)$ from $C^*(G_x,\Omega)$ to
$C^*(G,\Omega)$. Define the
 induced representation
of $(\tau,\rho_x)$ to be $L=(V,\rho_x)$, where
$V=\text{ind}_{G_x}^G(\tau)$, and $\rho_x$ is a representation of 
$C_0(\Omega)$ on
$H_V$, acting by $\rho_x(\phi)(g)(r)=\psi(r\cdot x)f(r)$.
We use the notation
$
L=(V,M)=\text{ind}_{(S,\Omega)}^{(G,\Omega)}(\tau,\rho_x).
$

\smallskip

We now present a different, equivalent way,
to induce representations of $C^*(G,\Omega)$.

\smallskip

Assume $H$ is a closed subgroup of $G$.

Assume that we have the $C^*$
algebra $C^*(G,\Omega)$,
and that $\pi$ is a representation of $C^*(H,\Omega)$
acting on the Hilbert
space $H_\pi$. Let $\xi$ and $\eta$ be arbitrary vectors in $H_\pi$.
Define the induced representation,
$L=\text{ind}_{(H,\Omega)}^{(G,\Omega)}(\pi)$ as follows:

 Let $B=C_c(H,\Omega)$, and define a $C_c(H,\Omega)$-valued
 inner product on 
 $C_c(G,\Omega)$ as follows:
$$
\langle f,g\rangle_{{}_B}(t,y)=
\int\limits_{s\in G}\overline{f}(s,s\cdot y)g(st,s\cdot y)
\dd\mu_{{}_G}(s)
$$
Define an inner product on 
$C_c(G,\Omega)\otimes H_{\pi}$ by

$$
\langle f\otimes\xi,g\otimes\eta\rangle_{{}_L}=
\bigl\langle\pi(\langle 
g,f\rangle_{{}_B})\xi,\eta\bigr\rangle{{}_{\tau}}
$$
\twoskip
Let $H_L$ be the completion of $C_c(G,\Omega)\otimes H_{\pi}$.

 The induced representation $L$ of $h\in C^*(G,\Omega)$ acts
 on the class of $f\otimes\xi$ 
by $(h*f)\otimes\xi$; see \cite{18}, page 204, or \cite{26}, 
page 340.

\twoskip
\proclaim{Proposition 1} Let $G$ be a group acting on a locally 
compact
Hausdorff space $\Omega$. Fix a point $x\in\Omega$.
Let $\tau$
a representation of $G_x$, and $\rho_x$ be a point evaluation of
$\Omega$. If $\tau$ is irreducible, then the induced
representation
$
L=\text{ind}_{(G_x,\Omega)}^{(G,\Omega)}(\tau,\rho_x),
$
of $C^*(G,\Omega)$,
is an irreducible representation of $C^*(G,\Omega)$.
\endproclaim

\demo{Proof}\enddemo
See \cite{26}, Proposition 4.2 on pg. 344.\qed

\twoskip
\noindent{\bf Definition 9.}
Let $L=(V,M)$ be a representation of $(H,\Omega)$, $s\in G$.
Let $K=sHs^{-1}$, and let $L^s=(V,M)^s=(V^s,M^s)$ be the covariant 
representation of
$(K,\Omega)$ given by
$
V^s(r)=V(s^{-1}rs),\text{and}  \ M^s(\phi)=M({}^s\phi).
$

\twoskip
The following has second-countability as a hypothesis, true in our case.

\twoskip
\proclaim{Theorem 2} If $S=\text{ind}_{(H,\Omega)}^{(G,\Omega)}(V,M)$
and $T=\text{ind}_{(K,\Omega)}^{(G,\Omega)}((W,N))$, then $S\cong T$
if and only if for some $s\in G$, we have
$T=\text{ind}_{(sKs^{-1},\Omega)}^{(G,\Omega)}((V,M)^s)$.
\endproclaim

\demo{Proof}\enddemo See \cite{15}, Theorem 2.1.\qed

\twoskip

We shortly prove an important proposition that enables us
to separate certain kernels of $C^*(G,\Omega)$. Several
definitions are needed.

\twoskip

The following equivalence relation $\sim_1$ was motivated by Dana P. Williams
\cite{26}, and reduces to his when $G$ is abelian.

\twoskip
\noindent{\bf Definition 10.}
Define an equivalence relation $\sim_1$ on
 $\go$ as follows: $(f,x)\sim_1 (f^\prime, x^\prime)$ when:

\itemitem{(1)}There exists $s\in G$ such that $x^\prime =s\cdot x$.

 \itemitem{(2)} For some $h\in \frak g_{x^\prime}^\perp$,
we have $\text{Ad}^*(s)f=f^\prime +h$.
\smallskip

We write $\Cal O_{(f,x)}$ for the equivalence class of the functional-point pair $(f,x)$.

Now we define a second quotient space $\sim$ on $\go$ by 
$\Cal O_{(f,x)}\sim\Cal O_{(h,y)}\ \iff\ 
\overline{\Cal O_{(f,x)}}=\overline{\Cal O_{(h,y)}}$. This is denoted
$\go/\sim$.

To verify
the above is an equivalence relation is a simple exercise.

\twoskip
Let $(V,M)$ be a covariant pair for $(G,\Omega)$. If $M$ is
given by $\rho_x$, a point evaluation of a point $x\in\Omega$, 
then for $\phi\in C_0(\Omega)$ we have
$\rho_x(\phi)(r)=\phi(r\cdot x)$.

\twoskip

Proposition 2 which follows Lemma 6 does 
not depend upon $G$ being a nilpotent Lie group. It has 
been proven by
Takesaki (\cite{25}, Theorem 7.2) when our
orbit space $\Omega/G$ has a $T_0$ topology. It was proven in
generality by Phil Green \cite{18}, pg. 210 Proposition 11, Part (ii). 
We present a different proof for our case.

\smallskip

As $\Omega$ is 
assumed separable, it is metrizable 
by \cite{20}, page 146 and Theorem 16, page 125.
Hence $\Omega$ is $T_4$, and the Tietze Extension Theorem may be 
applied to $\Omega$; this is needed to prove the next lemma.
\twoskip

\twoskip
\proclaim{Lemma 6} Let $x\in\Omega$, $G_x=$ the stabilizer of $x$ in $G$. Let
$C\subseteq G/G_x$ be compact in $G/G_x$. If $f\in C_c(G/G_x)$ is
supported on $C$, we may find a sequence of continuous functions 
$\{f_n\}_{n=1}^\infty\subseteq C_0(\Omega)$ with
$$
f_n(y)\to\left\lbrace
\aligned   f(s)\ \text{when}\ y=s\cdot x, \ \ y&\in C\cdot x \\
           0\ \ \ \ \ \ \ \ \ \ \ \ \ \ \ \ \ \ \ \ \ \ \ y&\notin C\cdot x. 
\endaligned \right.
$$

\endproclaim
\demo{Proof}\enddemo

We comment on the conclusion of this lemma. The function $h$ on $\Omega$ defined
by

$$
h(y)=\left\lbrace
\aligned  f(s)\ \text{when}\ y=s\cdot x, \ \ y&\in C\cdot x \\
           0\ \ \ \ \ \ \ \ \ \ \ \ \ \ \ \ \ \ \ \ \ \ \ \ y&\notin C\cdot x 
\endaligned\right.
$$
\twoskip
\noindent is not, in general, a continuous function on $\Omega$,
However, we may find a sequence of continuous functions limiting on $h$, so $h$
will be Borel.

The set $C\cdot x$ is closed and compact in $\Omega$ by sequential
arguments. Now simply define $f^\prime$ on $C\cdot x$ by
$f^\prime(s\cdot x)=f(s)$.
Assume that $\{y_n\}_{n=1}^\infty\subseteq C\cdot x$, and
$y_n\to y\in C\cdot x$. We need show that $f^\prime(y_n)\to f^\prime(y)$.
Assume 
that we have a collection $\bigl\{\{g_n\}_{n=1}^\infty,g\bigr\}\subseteq C$ and
$y_n=g_n\cdot x$ and $y=g\cdot x$, we may easily find such as $C$ is compact in
$G$ and $C\cdot x$ is compact in $\Omega$. As $C$ is compact, we may assume that
$g_n\to g^\prime$, and as $g^\prime\cdot x=g\cdot x$, and $\Omega$ is Hausdorff,
we have $g^\prime=g$. As $f$ is continuous on $C$, we have that $f^\prime$
is continuous
on $C\cdot x$ is clear by brute force. We now employ the Tietze Extension
Theorem to extend to all $\Omega$.

\twoskip
Now we note that $C\cdot x$ may not contain an open set. We find a sequence
$\{C_n\}_{n=1}^\infty$ of nested compact neighborhoods with
$\cap_{n=1}^\infty C_n=C\cdot x$. For each $n$, we use Urysohn's Lemma to find
a continuous function $g_n$ which is $1$ on $C\cdot x$ and $0$ on the closure of
$\widetilde{C_n\cdot x}$. We set $f_n=g_nf^\prime$ and we are done.\qed

\twoskip
\proclaim{Proposition 2} 
Let $x$ be any fixed point in $\Omega$, $\tau_1$ and $\tau_2$ be
two representations of $G_x$ such that $(\tau_1, \rho_x)$ and $(\tau_2, \rho_x)$
do not have the same kernel as representations of $C^*(G_x,\Omega)$.
 Defining $L_1=\text{ind}_{(G_x,\Omega)}^{(G,\Omega)}(\tau_1,\rho_x)$ and
$L_2=\text{ind}_{(G_x,\Omega)}^{(G,\Omega)}(\tau_2, \rho_x)$, 
we have $\text{ker}(L_1)\ne\text{ker}(L_2)$.
\endproclaim

\demo{Proof}\enddemo

We first define a new $C^*$ algebra, that being $C^*(G,G/G_x)$. We
observe that our orbit space $(G/G_x)/G$ for $C^*(G,G/G_x)$ {\it is}
a $T_0$
space; it consists of one point.  Hence $C^*(G,G/G_x)$
is Type I, and we have a canonical homeomorphism between
$\widehat{C^*(G,G/G_x)}$ and Prim$(C^*(G,G/G_x))$. 

We identify the point
$x\in\Omega$ with $x^\prime=e\in G/G_x$ in the new orbit space. We further 
observe that $G_{x^\prime}=G_x$ is clear. We use the same representations
$\tau_1$ and $\tau_2$ of $G_{x^\prime}$. Now define 
$$
\align
L_1^\prime&=\text{ind}_{(G_{x^\prime},\Omega)}^{(G,\Omega)}(\tau_1,\rho_{x^\prime})
\\ &\text{\ \ \ \ and} \\
L_2^\prime&=\text{ind}_{(G_{x^\prime},\Omega)}^{(G,\Omega)}
(\tau_2,\rho_{x^\prime}),
\endalign
$$ 
\noindent two irreducible
representations of $C^*(G,G/G_x)$, and that these have different kernels is
clear by the above comments and the fact that our orbit space is now $T_0$.

Throughout we assume that $\xi_{{}_i}$ and $\eta_{{}_i}$
$(i=1,2)$ are vectors in the space of $\tau_{{}_i}$ $(i=1,2)$. The 
representation
$(\tau_{{}_i},\rho_{x^\prime})$ of $C^*(G_{x^\prime},G/G_{x^\prime})$
we will refer to as $\pi_i$.

Now, as ker$(L_1^\prime)\ne$ker$(L_2^\prime)$, we may find an 
$h^\prime\in C^*(G,G/G_{x^\prime})$
with $h^\prime\in\text{ker}(L_1^\prime)$ and 
$h^\prime\notin\text{ker}(L_2^\prime)$.
So we may assume that for all elementary tensors 
$f^\prime\otimes\xi_{{}_1}$ that 
$L_1^\prime(h^\prime)(f^\prime\otimes\xi_{{}_1})=
h^\prime*f^\prime\otimes\xi_{{}_1}=0$
(or more properly, the class of
zero), and that for some elementary tensor $f^\prime\otimes\xi_{{}_2}$ that 
$L_2(h^\prime)(f^\prime\otimes\xi_{{}_2})\ne 0$. We may further assume
that $f^\prime$ is in the
dense subset of continuous functions with compact support having the form
$$
f^\prime(s,y)=\sum_{i=1}^l f^\prime_{1,i}(s)f^\prime_{2,i}(y).
$$ 
\noindent We may assume that each $f^\prime_{1,i}$ in the last displayed
formula is a continuous
function of compact support on $G$, and $f^\prime_{2,i}$ is a continuous 
function of
compact support on $G/G_x$ for each $i$, and $f^\prime$ and its component 
functions are
collectively supported on $C_1\times C_2$, where $C_1$ is a compact set in
$G$ and $C_2$ is a compact set in $G/G_x$. 

We will refer to all inner products of 
the $C^*$ algebra $C^*(G,G/G_x)$
as $\langle\cdot\,\ \cdot\rangle_{{}_i}^\prime$ to avoid confusion later on.

Re-iterating, we may assume by a
scaling argument that 
$$
\align
\langle h^\prime*f^\prime&\otimes\xi_{{}_1},f^\prime\otimes\xi_{{}_1}
\rangle_{{}_1}^\prime=0,\\
                          &\text{and} \\
\langle h^\prime*f^\prime&\otimes\xi_{{}_2},f^\prime\otimes\xi_{{}_2}
\rangle_{{}_2}^\prime=1,
\endalign
$$

Now assume that $\{h_n^\prime\}_{n=1}^\infty\subseteq C_c(G,G/G_x)$ and 
$h_n^\prime\to h^\prime$ in the topology of $C^*(G, G/G_x)$, and each 
$h_n^\prime$ has the 
form

$$
h^\prime_n(s,y)=\sum_{i=1}^{N_n}\phi_{n,i}^\prime(s)\psi_{n,i}^\prime(y),
$$
 
\noindent and each
is supported on a set $C_1^n\times C_2^n$, each $C_1^n$ a compact set in $G$ and
each $C_2^n$ a compact set in $G/G_x$.

We here note that $h_n^\prime*f^\prime$ (each $n$) is a continuous function
of compact support; this is proven on pages 32-33 of \cite{8}.

 We have 
$$
\aligned
\langle h^\prime_n*f^\prime\otimes\xi_{{}_1},f^\prime\otimes\xi_{{}_1}
\rangle_{{}_1}^\prime&\to 0 \\
\langle h^\prime_n*f^\prime\otimes\xi_{{}_2},f^\prime\otimes\xi_{{}_2}
\rangle_{{}_2}^\prime&\to 1,
\endaligned
$$
\twoskip
\noindent and more explicitly, when we \lq\lq untwist" the formulas we have

\twoskip

$$
\langle h^\prime_n*f^\prime\otimes\xi_{{}_i},f^\prime\otimes\xi_{{}_i}
\rangle_{{}_i}^\prime=
\bigl\langle\pi_i(\langle f^\prime,(h_n^\prime*f^\prime)
\rangle_{{}_{B_{{}_i}}}^\prime)\xi_{{}_i},
\xi_{{}_i}
\bigr\rangle_i^\prime=
$$
$$
\left\langle\pi_i\left(\left(\gamma_{{}_{G_x}}(t)\int\limits_{s\in G}
\overline{f^\prime}(s,s\cdot y)(h^\prime_n*f^\prime)(st,s\cdot y)
\dd\mu_{{}_G}(s)\right)(t,y)\right)\xi_{{}_i},\xi_{{}_i}\right\rangle_i^\prime=
$$
$$
\multline 
\int\limits_{r\in G_x} \biggl(\gamma_{{}_{G_x}}(r)
\int\limits_{s\in G}\overline{f^\prime}(s,s\cdot x^\prime) \\
(h_n^\prime *f^\prime)(sr,s\cdot x^\prime)\dd\mu_{{}_{G}}(s)\biggr)
\bigl\langle \tau_i(r)\xi_{{}_i},\xi_{{}_i}\bigr\rangle_{{}_i}^\prime
\dd\mu_{{}_{G_x}}(r)\underset{n\to\infty}\to{\longrightarrow}
\endmultline 
$$
$$
\to\left\lbrace\aligned 
 0\ \ \ \ i&=1 \\ 
 1\ \ \ \ i&=2.
\endaligned \right.
$$

\twoskip     

Now we return to $C^*(G,\Omega)$. All inner products in this $C^*$ algebra we
refer to by $\langle\cdot\,\ \cdot\rangle_{{}_i}$ (no ${}^\prime$'s
so as to distinguish from those in $C^*(G,G/G_x)$\ ). 

Again, we have the hypothesis that
$\tau_1\not\cong\tau_2$ as representations of $G_x$. 
Let 
$$
\align 
L_1&=\text{ind}_{(G_x,\Omega)}^{(G,\Omega)}(\tau_1, \rho_x),\\
   & \text{ \ \ \ and} \\
L_2&=\text{ind}_{(G_x,\Omega)}^{(G,\Omega)}(\tau_2, \rho_x).
\endalign  
$$

\noindent We define a collection
of functions $\bigl\{\{h_{n,j}\}_{n,j=1}^\infty, \{f_j\}_{j=1}^\infty\bigr\}$ 
in $C_c(G,\Omega)$ 
in a special way. Define

$$
f_j(s,y)=\sum_{i=1}^lf_{1,i}(s)f_{2,i,j}(y),
$$

\noindent
each $f_{1,i}=f^\prime_{1,i}$; we do not change
these functions. Now we choose $f_{2,i,j}$ to satisfy 
$f_{2,i,j}(s\cdot x)=f^\prime_{2,i}(s\cdot x^\prime)$, so $f_{2,i,j}$
\lq\lq behaves the same as $f^\prime_{2,i}$" on the specific set $C_2\cdot x$,
extend to all $\Omega$, and require 

$$
f_{2,i,j}(y)\underset{j\to\infty}\to{\longrightarrow}\left\lbrace
\aligned 
f^\prime_{2,i}(s\cdot x^\prime)\ &\text{when}\ y=s\cdot x,
\ \ s\in C_2 \\
0\ \ \ \ \ \
\   \ &\text{otherwise},
\endaligned \right.
$$

\noindent see
Lemma 6 for specifics. We note that we have
 
$$
\lim_{j\to\infty}f_j(s,y)=\left\lbrace
\aligned
f^\prime(s,s\cdot x^\prime)\ &\text{when}\ y=s\cdot x,
\ \ s\in C_2 \\
0\ \ \ \ \ \
\   \ &\text{otherwise},
\endaligned\right.
$$

\noindent Similarly use Lemma 6 to define a sequence
$\{h_{n,j}\}_{n,j=1}^\infty$ with (for each fixed $n$)
$$
h_{n,j}(s,y)\underset{j\to\infty}\to{\longrightarrow}\left\lbrace
\aligned 
h^\prime_n(s,s\cdot x^\prime)\ &\text{when}\ y=s\cdot x,
\ \ s\in C_2 \\
0\ \ \ \ \ \
\ \ &\text{otherwise}.
\endaligned \right.
$$
For each $n$ and all
$j$ (fixed $n$) we may assume
that the functions $f_j$ and $h_{n,j}$ have support contained in the set
$S^n_1\times S^n_2$, $S^n_1$ compact in $G$, and $S^n_2$ compact in $\Omega$;
see again the proof of Lemma 6. We also
(for each fixed $n$) choose them uniformly bounded for all $j$.

\twoskip

Now we re-write our old calculations with these new functions and follow
the same steps now in the $C^*$ algebra $C^*(G,\Omega)$:

$$
\langle h_{n,j}*f_j\otimes\xi_{{}_i},f_j\otimes\xi_{{}_i}\rangle_{{}_i}=
$$
$$
\int\limits_{r\in G_x} \left(\gamma(r)
\int\limits_{s\in G}\overline{f_j}(s,s\cdot x)
(h_{n,j}*f_j)(sr,s\cdot x)\dd\mu_{{}_{G}}(s)\right)
\bigl\langle \tau_i(r)\xi_{{}_i},\xi_{{}_i}\bigr\rangle_{{}_i}
\dd\mu_{{}_{G_x}}(r).\tag1
$$

Now we observe the integral in $s$ is over a compact set as 
$f_j$ is of compact support, and by the way that we choose the functions
$h_{n,j}*f_j$ for all $j$ these have support
contained in another compact set we denote $K_1\times K_2$. 
So we must have $sr\in K_1$, and as $s\in S_1^n$ is forced already, we have
$r\in ((S_1^n)^{-1}\cdot K_1)\cap G_x$ is forced, and this is compact in $G_x$.

Now we observe we may apply the Bounded Convergence Theorem to
the above integral, these functions are converging in $j$ and for each $n$
have been chosen in a bounded fashion with compact support, hence are easily
bounded above by an integrable function. So for each fixed $n$, 
as we limit $j\to\infty$, formula (1) above converges to
$$
\left\lbrace\aligned 
\langle h^\prime_n*f^\prime\otimes\xi_{{}_1},f^\prime\otimes\xi_{{}_1}
\rangle_{{}_1}^\prime\
\ \ (i&=1) \\ 
\langle h^\prime_n*f^\prime\otimes\xi_{{}_2},f^\prime\otimes\xi_{{}_2}
\rangle_{{}_2}^\prime \ 
\ \ (i&=2),
\endaligned \right.
$$
\noindent where we remind the reader that these inner products are in our
$C^*$ algebra $C^*(G,G/G_x)$.
\noindent Furthermore, we may assume (fixed $n$) that for all $j\ge n$
$$
\aligned 
\bigl|\langle h_{n,j}*f_j\otimes\xi_{{}_1},f_j\otimes\xi_{{}_1}\rangle_{{}_1}-
\langle h^\prime_n*f^\prime\otimes\xi_{{}_1},f^\prime\otimes\xi_{{}_1}
\rangle_{{}_1}^\prime\bigr|&<\frac1n 
 \\ 
\bigl|\langle h_{n,j}*f_j\otimes\xi_{{}_2},f_j\otimes\xi_{{}_2}\rangle_{{}_2}-
\langle h^\prime_n*f^\prime\otimes\xi_{{}_2},f^\prime\otimes\xi_{{}_2}
\rangle_{{}_2}^\prime\bigr|&<\frac1n,
\endaligned
$$

\noindent where we have mixed inner products of $C^*(G,\Omega)$ and
$C^*(G,G/G_x)$ in the above.

Now we may choose the diagonal sequence $\{h_{j,j}\}_{j=1}^\infty$ and the
sequence $\{f_j\}_{j=1}^\infty$ in the above integrals, and for $i=1$ the
sequence converges to 0, for $i=2$ it converges to $1$. As $\xi_1$ was 
arbitrary, we now have
$L_1(h_{j,j})\to 0$ and $L_2(h_{j,j})\not\to 0$, showing
that $L_1$ and $L_2$ cannot have the same kernel, the desired result.\qed

\twoskip
\noindent{\bf Corollary 1.}
Assume 
the pairs $(f_1,x)$ and $(f_2, x)$ give rise to
irreducible representations 
$(\tau_1, \rho_x)$ and $(\tau_2,\rho_x)$ of
$C^*(G_x,\Omega)$. These induce to equivalent irreducible 
representations
of $C^*(G,\Omega)$ if and only if 
$\tau_1$
and $\tau_2$ are the equivalent irreducible representations of $G_x$.
This says that
$(f_1,x)$ and $(f_2,x)$ are
in the same equivalence class mod $\sim_1$, and $f_1|_{\frak g_x}$ and 
$f_2|_{\frak g_x}$
are in the same $G_x$ orbit in $\frak g_x^*$; see Definition 10.

\twoskip

For the remainder of this section, we use Siegfried Echterhoff 
\cite{7} as our main reference.

\twoskip
\noindent{\bf Definition 11.}
Let $G$ be our connected, simply-connected nilpotent Lie group, and
let $\Cal K(G)$ be the space of closed subgroups of $G$.
 Let $\Cal N$ be a locally compact space. Assume 
$H:\Cal N\mapsto\Cal K(G)$ and $H(i)= H_i$ is a  continuous map from
$\Cal N$ to $\Cal K(G)$. Define:

$$
\Cal N^H=\{(i,x)\in\Cal N\times G\ |\ x\in H_i\} 
$$

\twoskip

In this paper, we use
$\Cal N=\Bbb N\cup\infty$, the one-point compactification of the 
natural numbers.

\twoskip

The following results are from \cite{7}.

\twoskip
\noindent{\bf Definition 12.}

Now let $G$ again be a nilpotent Lie group and $(G,\Omega)$ denote a covariant 
system. We make the space $C_c(\Cal N^H, C_0(\Omega))$ 
into a normed ${}^*$- algebra.
Define multiplication, involution, and norms by

$$
\align
f*g(i,t,x) &=\int\limits_{s\in H_i}f(i,s,x) g(i,s^{-1}t, s^{-1}\cdot x)
\dd\mu_{{}_{H_i}}(s) \\
f^*(i,t,x) &= \overline{f}(i,t^{-1},t^{-1}\cdot x) \\
\Vert f\Vert_1 &=\underset{{}_{i\in\Cal N}}\to{\text{sup}}\
\int\limits_{s\in H_i}\underset{x\in\Omega}\to{\text{sup}}\ |f(i,s,x)|
\dd\mu_{{}_{H_i}}(s)
\endalign
$$

\noindent for $f,g\in C_c(\Cal N^H,C_0(\Omega))$

\twoskip
\noindent{\bf Definition 13.}
Denote by $L^1(\Cal N^H, C_0(\Omega))$
\comment
 and $L^1(\Cal N^K_H,C_0(\Omega))$ 
\endcomment
the
completion of $C_c(\Cal N^H,C_0(\Omega))$
\comment
and 
$C_c(\Cal N^K_H, C_0(\Omega))$
\endcomment
 with respect to the above norm.
Any covariant representation $\tau$ of
the system $\bigl(H_i,C_0(K_i,C_0(\Omega))\bigr)$ defines a
${}^*$ - representation of $L^1(\Cal N^H,C_0(\Omega))$ by

$$
(i,\tau)(F) = \tau(F_i)
$$

\twoskip
\noindent{\bf Definition 14.}
The representation space

$$
\Cal R(C^*(\Cal N^H, C_0(\Omega)))=\bigl\{(i,\rho)\ |\ i\in\Cal N,
\rho\in\text{Rep}(H_i,C_0(\Omega))\bigr\}
$$

\noindent  with the relative topology of
$\text{Rep}(C^*(\Cal N^K,C_0(\Omega)))$ 
 is the {\it subgroup representation
space} of
$C^*(\Cal N^K,C_0(\Omega))$.

Let $C^*(\Cal N^H, C_0(\Omega))$ be the enveloping $C^*$ algebra
of $L^1(\Cal N^H,C_0(\Omega))$
 We call 
$C^*(\Cal N^H, C_0(\Omega))$
 the {\it subgroup algebra} of
$(\Cal N^H,C_0(\Omega))$.
\twoskip

A standard technique in nilpotent harmonic analysis is induction
from codimension one subgroups. The next proposition extends this idea
from the above-defined algebras to our transformation
group $C^*$ algebras.

\twoskip

\proclaim{Proposition 3} Let $G$ be a nilpotent Lie group acting upon
the locally compact Hausdorff space $\Omega$.

\item{(1)}{Assume we have two continuous maps,
$k(n)$ and $h(n)$,
 from $\Cal N=\Bbb N\cup\infty$ to $\Cal K(G)$.
Also assume that for all $n$, 
$h(n)=H_n$ is codimension $1$ in $k(n)=K_n$.} 

\item{(2)}{Assume we have a collection of representations
$\{\{(i,\pi_i)\}\ |\ i\in\Cal N\}$ in the space
 $\Cal R(C^*(\Cal N^K,C_0(\Omega)))$, and,
for each $i\in\Bbb N$, 
 $\pi_i$ is induced from $(H_i,\Omega)$. Assume for each $i\in\Bbb N$
that $\pi_i$ is
irreducible on $C^*(K_i,\Omega)$, and that each $\pi_i$ corresponds
to a functional-point pair $(f_i,x_i)$.}

\item{(3)}{Assume in the space
 $\Cal R(C^*(\Cal N^K,C_0(\Omega)))$ that
$(i,\pi_i)\to(\infty,\pi_\infty)$.}

\twoskip

Then,

\twoskip

\item{(1)}{By passing to a sub-sequence,
we may choose a collection
$
\left\{\sigma_{i}\right\}_{i=1}^\infty,\
$
with each $\sigma_{i}$
an irreducible representation of $C^*(H_{i}, \Omega)$, and
$\sigma_\infty$, an
irreducible
representation of $(H_\infty,\Omega)$,
these representations
satisfying:

$$ 
\text{ker}(\pi_{i})=
\text{ker}(\text{ind}_{(H_{i},\Omega)}^{(K_{i},\Omega)}
(\sigma_{i})),
$$

$$
\pi_\infty\prec\text{Ind}_{(H_\infty,\Omega)}^{(K_\infty,\Omega)}
(\sigma_\infty),\ \text{and}\
(i,\sigma_{i})\to(\infty,\sigma_\infty)
\ \text{in}\ \Cal R(C^*(\Cal N^H,C_0(\Omega))).
$$}

\item{(2)}{For each $i$, the functional-point
pairs (Definition 10) corresponding to $\sigma_{i}$ may be chosen
in the same 
 equivalence class mod $\sim_1$ as $(f_i,x_i)$.}

\endproclaim

\demo{Proof}\enddemo By assumption, for each
$i\in\Bbb N$, $\pi_i$ is induced from
$(H_i,\Omega)$, so assume
that $\pi_i=\text{ind}_{(H_i,\Omega)}^{(K_i,\Omega)}(\sigma_i)$.
By Corollary 1, we may assume that the functional-point pairs 
(Definition 10) 
corresponding to $\sigma_i$ are the same as those corresponding to $\pi_i$.

The restriction map
$\Cal R(C^*(\Cal N^K,C_0(\Omega)))\mapsto\Cal R(C^*(\Cal N^H,C_0(\Omega)))$
is continuous (Proposition 7 of \cite{7}),
so we have:

$$
(i,\pi|_{(H_i,\Omega)})\underset{i\to\infty}\to{\longrightarrow}
 (\infty,\pi|_{(H_\infty,\Omega)}).
$$

Note by \cite{3}, Lemma 1.1.8 that $H(i)$ is normal in $K(i)$. By
\cite{26}, Proposition 7 on pg. 70, for each $i$, the
restricted representation is equivalent to the following direct 
integral:

$$
(i,\pi_i|_{(H_i,\Omega)})
=(i, \int\limits^\oplus_{{}_{s\in 
K_i/H_i}}
\!\sigma_i^s \ \ ).
$$

Prior to finishing the proof, we present several facts particular to connected, simply connected 
nilpotent Lie groups. We observe our group representations are 
not generally irreducible, but are induced from irreducible representations
of subgroups, hence are direct integrals of irreducibles, as are their
restrictions to subgroups.

\smallskip

{\bf Fact 1}: Let $\pi=\pi_f$ be an infinite-dimensional 
irreducible representation of $G$. Let $G_0$ be 
a codimension one subgroup of $G$, and let $f_0=f|_{\frak g_0}$ the 
restriction of
$f$ to $\frak g_0$. Let $\pi_0=\pi_{f_0}$ be the irreducible representation
of $G_0$ associated by Kirillov theory to $f_0$ on $\frak g_0$.

By \cite{3}, Theorem 2.5.3,
the representation $\pi|_{G_0}$ is either a representation that
induces directly to $\pi$, or is the following direct integral:
$$
\pi|G_0\ = \int\limits^\oplus_{{}_{s\in G/H}}
\! (\pi_0)^s,
$$
\noindent and $\pi = \text{ind}_{G_0}^G(\pi_0)$.

Let $\sigma = \text{ind}_{G_0}^G(\pi_0)$. We have one of the following:

\noindent {\it a)} $\sigma\cong\pi$

\noindent {\it b)} $\sigma\cong\infty\cdot\pi$

{\bf Fact 2}: Recall Lemma 2.4 on pg. 338 of \cite{26}:

\smallskip

{\bf Lemma.} {\it Let $A$ be a $C^*$ algebra. 
Let $\{I_\alpha\}_{\alpha\in\Lambda}$ be a net of ideals
of $A$, converging to $I$. 
Assume that each $I_\alpha$ is the intersection
of a set of primitive ideals
$F_\alpha\in\Cal K(${\it Prim}$(A))$, and $F$ corresponds to $I$.
Then, given any $F^\prime\supseteq F$, there is a
subnet of primitive ideals of $A$, 
$\{I_\beta\}_{\beta\in\Lambda^\prime}$, such that there are
$F^\prime_\beta\in F_\beta$, with $\{F^\prime_\beta\}_{\beta\in\Lambda^\prime}$ converging
to $F^\prime$ in Prim$(A)$.}

\smallskip
\twoskip

Now we return to the proof of our proposition.

\smallskip

Note that $\pi_\infty$ may be
one-dimensional, hence not induced from $H_\infty$.

By the first of the above two facts, choose an irreducible
representation $\sigma^\prime_\infty$ of $C^*(H_\infty,\Omega)$ with
$\sigma^\prime_\infty\prec\pi_\infty|_{(H_\infty,\Omega)}$ and
$\pi_\infty\prec\text{ind}_{(H_n,\Omega)}^{(K_n,\Omega)}
(\sigma^\prime_\infty)$. 

We may assume
that the functional-point pair for $\sigma^\prime_\infty$ is in
the same $\sim_1$ class as that of $\pi_\infty$.

By an application of the second fact above, 
we may pass
 to a sub-sequence, and choose
$\{(n, \sigma_{n}^\prime)\}_{n=1}^\infty$ with
$(n,\sigma_{n}^\prime)\to(n,\sigma^\prime_\infty)$.

By an application of Theorem 2 and Corollary 1, for all $n\in\Bbb N$, we have
$\sigma_{n}^\prime\cong\sigma^{s_{n}}_{n}$
for some $s_{n}\in H_{n}$. So the 
functional-point pairs associated to $\sigma_{n}$ and
$\sigma_{n}^\prime$ are in the same equivalence class
mod $\sim_1$. \qed

\twoskip
\proclaim{Lemma 7} Assume that $G$ is a nilpotent Lie group acting on
a locally compact Hausdorff space $\Omega$.

Assume we have a
continuous map $h(i)=H_i$ from
$\Cal N=\Bbb N\cup\infty$ to $\Cal K(G)$, and
$H_i\to H_\infty$.

Assume we have a collection of one-dimensional representations 
$\bigl\{\{\pi_i\ |\ i\in\Cal N\}\bigr\},\ \pi_i\in\text{Rep}(C^*(H_i,\Omega))$,
each $\pi_i=(\chi_{f_i}, \rho_{x_i})$.

Assume in the space
$\Cal R(C^*(\Cal N^H,C_0(\Omega)))$ that
$(i,\pi_i)\to(\infty,\pi_\infty)$.

\twoskip

Then:

\smallskip

\item{(a)} $x_n \to x$, and

\smallskip

\item{(b)} by passing to a sub-sequence, we may choose another
collection $\{f^\prime_i\}_{i=1}^\infty$ 
with
$
\chi_{f^\prime_i} = \chi_{f}.
$

\smallskip

Furthermore, for each $i\in\Cal N$, $(f_i^\prime, x_i)$
may be chosen in the
same $\sim_1$ class of $\frak g^*\times\Omega$ as $(f_i, x_i)$.

\endproclaim

\demo{Proof}\enddemo
For any $\phi\in C_c(G)$ and $\psi\in C_c(\Omega)$, we may consider
$\phi\cdot\psi\in C^*(\Cal N^H,C_0(\Omega))$ by setting
$(\phi\cdot\psi)(i)=\phi|_{H_i}\cdot\psi$.
By hypothesis,
$$
\multline
\pi_i(\phi|_{H_i}\cdot\psi)=
\psi(x_i)\int\limits_{s\in H_i}\phi(s)\chi_{f_i}(s)
\dd\mu_{H_i}(s) \\
\longrightarrow
\psi(x)\int\limits_{s\in H}\phi(s)\chi_{f}(s)\dd\mu_H(s)=
\pi_\infty(\phi|_{H_\infty}\cdot\psi).
\endmultline
$$

\noindent Let $e$ denote the identity element of the group. By
the continuity of
$\text{Res}_e^{H_n}$, the
the restriction map from $\Cal R(C^*(\Cal N,C_0(\Omega)))$ to 
$\Omega$,
  (Proposition 7 of \cite{7}),
 we have
for $\phi\in C_0(\Omega)$,
$\text{Res}_e^{H_n}(\pi_i)(\phi)=\phi(x_i)$, and $\phi(x_i)\to\phi(x)$
for all $\phi\in C_0(\Omega)$; consequently
$x_i\to x$.

Choose $\psi$ identical one in a neighborhood of 
$x\in\Omega$.
  In Fell's
subgroup-pair topology \cite{12} we have:

$$
\langle \chi_{f_i},H_i\rangle\to
\langle\chi_{f},H_\infty\rangle.
$$

\noindent Pass to a sub-sequence, and 
use Lemma 5 
to find a sequence
$\{f^\prime_i\}_{i=1}^\infty\subseteq\frak g^*$ with
$\chi_{(f^\prime_i|_{H_i})}=\chi_{(f_i|_{H_i})}$ and $f_i^\prime\to 
f_\infty$ in
$\frak g^*$ and 
 by Lemma 5, each $f^\prime_i$ is in the same $H_i$-orbit as
$f_i$ and $(f_i^\prime, x_i)$ is in the same $\sim_1$ equivalence class
of $\frak g^*\times\Omega$.\qed

\twoskip

Now for our \lq\lq big" lemma, important in the next section.

\twoskip
\proclaim{Lemma 8}
Assume that we have a collection 
$\bigl\{\{F_i,\}_{i=1}^\infty,F\bigr\}$ of
primitive ideals of $C^*(G,\Omega)$, with
$$\multline
F_i=\text{ker}(L_i)=
\text{ker(ind}_{(G_{x_i},\Omega)}^{(G,\Omega)}
(\tau_{{}_{f_i,x_i}}, \rho_{x_i})) \\
\longrightarrow F=\text{ker}(L)=
\text{ker(ind}_{(G_x,\Omega)}^{(G,\Omega)}(\tau_{{}_{f,x}},\rho_x)).
\endmultline
$$

\noindent By passing to a sub-sequence we may choose:

$$
\{y_{i}\}_{i=1}^\infty\subseteq\Omega  \
\text{and}\ \{f^\prime_{i}\}_{i=1}^\infty \subseteq\frak g^*,\
\text{with}\ y_{i}\to y\in\Omega, f^\prime_{i}\to 
f^\prime,
$$

\noindent with

$$
\text{ker}(L_{i})=
\text{ker(ind}^{(G,\Omega)}_{(G_{y_{i}},\Omega)}
(\tau_{{}_{f^\prime_{ii},y_{i}}}, \rho_{y_{i}}))
$$

\noindent and

$$
\text{ker}(L)=\text{ker}(\text{ind}_{(G_{x^\prime},\Omega)}^{(G,\Omega)}
(\tau_{f^\prime}, \rho_{y})).
$$

Furthermore, the functional-point pairs $(f^\prime_i,y_i)$
corresponding to $\tau_{i}^\prime$ 
may be chosen in the same $\sim$
equivalence class as $(f_i,x_i)$, and $(f^\prime,y)$ in the same
$\sim$ equivalence class as $(f,x)$.

\endproclaim

\demo{Proof}\enddemo
By Lemma 4.5 of \cite{26} we may
assume that $x_i\to x$ and that the
$\sim$ equivalence classed don't change.  

If an infinite sub-sequence of $\{L_i\}_{i=1}^\infty$
consists of one-dimensional representations of $C^*(G, \Omega)$,
we may use the constant sequence $G_{x_i} = G$ for all $i$
and Lemma 7.

 If no infinite
sub-sequence of $\{L_i\}_{i=1}^\infty$ 
consists of one-dimensional representations,
we may assume that each is induced from a codimension
one subgroup $H_i$.
We pass to a sub-sequence and assume that for the sequence $\{\frak p_i\}_{i=1}^{\infty}$ of
polarizing subalgebras, we have dim$(\frak p_i)$ constant.

 By compactness of $\Cal K(G)$, we
may pass to a sub-sequence and assume that
$G_{x_i}\to S\subseteq G_x$ and $\frak p_i\to \frak p$, where $\frak p$ may not be
polarizing for the action of $f$ on $\frak g_x$.
Assume that 
a representation at least weakly containing $L$ is
induced from $(H,\Omega)$.

We may use Proposition 3 to successively reduce our problem by one
dimension, until we do have a sequence of one-dimensional
representations and Lemma 7 may be applied to these.\qed

\vfill\eject

\topmatter
\title\chapter{2} \endtitle
\endtopmatter

\title The Topology on Prim$(C^*(G,\Omega))$\endtitle
\endtopmatter

\twoskip
\proclaim{Lemma 9} Let $x\in\Omega$, $f\in\frak g^*$ be given.

Let $\frak p_x$ be polarizing for the restriction of $f$
to $\frak g_x$, and let $\frak p$ be isotropic (not necessarily 
polarizing)
for the restriction of $f$ to $\frak g_x$.  Then
$L^\prime=(V^\prime,\rho_x)= \text{ind}_{(P,\Omega)}^{(G,\Omega)}
(\chi_{{}_{f,P}},\rho_x)$ weakly contains $L=
(V,M)=\text{ind}_{(P_x,\Omega)}^{(G,\Omega)}(\chi_{{}_{f,P_x}},\rho_x)$.
\endproclaim

\demo{Proof}\enddemo

Let $\chi_{f,P}$ be the character of $P$ determined by $f$,
and let
$\chi_{f,P_x}$ be the corresponding character of $P_x$.
By Lemma 4, on the level of stabilizer subgroups, we
have 
$
\text{ind}_{(P_x,\Omega)}^{(G_x,\Omega)}(\chi_{f,P_x},\rho_x)\prec
\text{ind}_{(P,\Omega)}^{(G_x,\Omega)}(\chi_{f,P},\rho_x).
$
As induction preserves weak containment
(Proposition 9, \cite{18}), the conclusion is clear from  
this and
\lq\lq induction in stages"; see Proposition 8, pg. 207 of 
\cite{18}.\qed

\twoskip
\noindent {\bf Definition 14.}
      Define $\phi:\go\mapsto\pgo$ by
$$
\phi 
(f,x)=\text{ker(ind}^{(G,\Omega)}_{(G_x,\Omega)}(\tau_{f,x},\rho_x)).
$$
\twoskip
\proclaim{Lemma 10} $\phi$ is continuous in the product topology of
$\frak g^*\times\Omega$.\endproclaim

\demo{Proof }\enddemo

Assume in the product topology of $\frak g^*\times\Omega$ that
$(f_n,x_n)\to(f,x)$. Denote the sequence of stability subgroups
as $\left\{G_{x_n}\right\}_{n=1}^\infty$; we assume that
$G_{x_n}\to S\subseteq G_x$. 

Let $\left\{\frak p_n\right\}_{n=1}^\infty$ denote the sequence 
of polarizing
subalgebras of $\frak g_{x_n}$ for the sequence of functionals
$\left\{f_n\right\}_{n=1}^\infty$. As $f_n\to f$, we may assume 
that $\frak p_n\to\frak p$, where $\frak p$ may be of lower
dimension than a polarizing
subalgebra $\frak p_x$ of $\frak g_x$ for $f|_{g_x}$.

Denote by $g$ the restricted
functional $f|_{\frak p}$. Denote by $\sigma$
the representation of $C^*(P,\Omega)$ given by the obvious
character $\chi_g$
of $P$ and a point evaluation $\rho_x$ of $\Omega$.

For each $n$ denote by $\pi_n$ the representation of $C^*(P_{x_n},\Omega)$ given by
the pair $(\chi_{f_n}, \rho_{x_n})$, and by $\pi$ the representation of
$C^*(P_x,\Omega)$ given by the pair $(\chi_f,\rho_x)$. 

Define $L_n=\text{ind}_{(G_{x_n},\Omega)}^{(G,\Omega)}(\pi_n)$,
an irreducible representation of $C^*(G,\Omega)$, and
$L=\text{ind}_{(P_{x},\Omega)}^{(G,\Omega)}(\rho_x)$,
also an irreducible representation of $C^*(G,\Omega)$. 

Let
$L^\prime=\text{ind}_{(P,\Omega)}^{(G,\Omega)}(\sigma)$ be
induced from $C^*(P,\Omega)$; this may not be irreducible.

Amending Echterhoff \cite{7}, Proposition 6 on pg. 69 slightly for
our purposes, we have:

{ \it Let $(G,\Omega)$ be a covariant system, $\Cal N$ be the locally compact
space $\Bbb N\cup\infty$, and $P(n)=\frak p_{x_n}$ be a continuous map
with $P(n)=\frak p_{n}\to \frak p=P(\infty)$ in $\Cal K(G)$. Then the map

$$
Ind_{P}^G:\Cal R(C^*(\Cal N^p,\Omega))\mapsto \widehat{C^*(G,\Omega)};
\ \ (n, \pi_n)\mapsto 
(n, \text{ind}_{(P_n,\Omega)}^{(G,\Omega)}(\pi_n))
$$
\noindent is continuous.  }

Using this, we have $L_n\to L^\prime$, but as $L\prec L^\prime$
(Lemma 9), we are done.\qed

\twoskip
\noindent{\bf Definition 15.}
We say that $C^*(G,\Omega)$ is {\it EH regular} if:

\item{(1)}{$C^*(G,\Omega)$ is quasi-regular,}
\item{(2)}{for every $P\in$ Prim$(C^*(G,\Omega))$, there is an
$x\in\Omega$ and an irreducible representation $\tau$ of $G_x$
such that  $P=\text{ker}(\text{ind}_{(G_x,\Omega)}^{(G,\Omega)}
(\tau,\rho_x))$.}

This has been established in our case by Jon Rosenberg and Elliot 
Gootman \cite{17}.

\twoskip

Remember our equivalence relation $\go / \sim$, as well as the
equivalence class $\Cal O_{(f,x)}$ of $(f,x)$ in $\go$;
 see Definition 10.

\twoskip
\noindent{\bf Definition 16.}
Using $\sim$ from Definition 10, define $\psi\ :\go/\sim\ \mapsto\text{Prim}(C^*(G,\Omega))$ by
$\psi(f,x)=\text{ker}(\text{ind}_{(G_x,\Omega)}^{(G,\Omega)}
(\tau_{f,x},\rho_x)$. We show that $\psi$ factors through $\sim$
in the next lemma.

\twoskip
\proclaim{Lemma 11} The map $\psi$ is one-to-one onto 
on $\go/\sim$.
\endproclaim
\demo{Proof}\enddemo

Onto is clear by EH regularity; we show that $\psi$ is 1-1.

When $\overline{\Cal O_{(f,x)}}=
\overline{\Cal O_{(h,y)}}$, that
the primitive ideals defined by $(f,x)$ and $(h,y)$ are the same
follows from the continuity of the map $\phi$ (Lemma 10)
and the fact that primitive ideal spaces are $T_0$.

Now assume that $\psi(f,x)=\psi(h,y)$.

As $\overline{G\cdot x}=\overline{G\cdot y}$, by
Lemma
4.5 \cite{26}, we may find a sequence
$\{g_n\}_{n=1}^\infty\subseteq G$ such that $g_n\cdot x\to y$.

Define $f_n=\text{Ad}^*(g_n)f$, and $x_n=g_xn$.

We have: 
$
\psi(f_n,x_n)=\psi(f,x),
$
and the sequence $\psi(f_n,x_n)$ always remains
in the $\sim_1$ equivalence class  of $(f,x)$, and 
$\psi(f_n,x_n)=\psi(h,y)$.
\twoskip

If a sub-sequence of the sequence 
$\psi((f_n,x_n))$ consists
of kernels of one-dimensional representations of $C^*(G,\Omega)$, 
we may use Lemma 7 to get an 
equivalent sequence $(f_n, x_n)$ converging
to $(h,y)$ with each $(f_n, x_n)$ still in 
the equivalence class of
$(f,x)$.

\smallskip

Thus we may assume there is no 
sub-sequence of one-dimensional representations.

Note that
$\left\{\psi(f_n,x_n)\right\}_{n=1}^\infty$ is a
constant sequence of ideals equal to $\psi(h,y)$.

In Lemma 8  we showed that if we
 had a collection of primitive
ideals, 
$\bigl\{\{F_n,\}_{n=1}^\infty, F\bigr\}$ of
$C^*(G,\Omega)$, with
$$
F_n=
\text{ker(ind}_{(G_{x_n},\Omega)}^{(G,\Omega)}
(\tau_{{}_{f_n,x_n}}, \rho_{x_n})) 
\longrightarrow F=
\text{ker(ind}_{(G_x,\Omega)}^{(G,\Omega)}(\tau_{{}_{f,x}},\rho_x)),
$$

\noindent that, passing to a sub-sequence, we could choose:

$$
\{y_{n}\}_{n=1}^\infty\subseteq\Omega  \
\text{and}\ \{f^\prime_{n}\}_{n=1}^\infty \subseteq\frak g^*,\
\text{ with}\ y_{n}\to y\in\Omega, f^\prime_{n}\to 
f^\prime,
$$

\noindent with
$
F_n=
\text{ker(ind}^{(G,\Omega)}_{(G_{y_{n}},\Omega)}
(\tau_{{}_{f^\prime_{n},y_{n}}}, \rho_{y_{n}}))
$
 and
$
F=\text{ker}(\text{ind}_{(G_{y},\Omega)}^{(G,\Omega)}
(\tau_{f_{f, y}^\prime}, \rho_{y})).
$
We may also choose the functional-point
pairs corresponding to $\tau_{f^{\prime}_n}$ and $\tau_{f^\prime}$ 
in the same $\sim_1$
equivalence class as those corresponding to $\pi_{n}$.

Thus by choosing appropriate functional-point pairs
in the equivalence class of $(f,x)$ we may realize 
$(h,y)$ as a subsequential limit, and
the equivalence class closures are equal.\qed

\twoskip
Now to characterize the primitive ideal space of $C^*(G,\Omega)$.

\twoskip
\proclaim{Theorem 3} The map $\psi$ is a homeomorphism from
$\go/\sim$ to \linebreak Prim$(C^*(G,\Omega))$.
\endproclaim

\demo{Proof }\enddemo 
We follow the philosophy of Dana P. Williams' Theorem 5.3
of \cite{26}.

We have the diagram:

$$
\CD
\go            @.          {}  \\
@VVqV                      @.  \\
\go /\! \sim    @>\psi >>      \pgo
\endCD
$$

As $\phi$ (Lemma 10) and the natural map $q$  of $\go$ to
$\go/ \sim$ are continuous, $\psi$ is 
continuous. By the last lemma $\psi$ is 1-1 onto.

Let $F$ be closed in $\go$ and saturated with
respect to $\sim$. We show that $\psi (F)$ is closed in 
$\pgo$.

   Assume $\{A_n\}_{n=1}^\infty\subseteq\psi (A)$ and $A_n\to 
A$.
   By  EH regularity, assume that $A_n=\text{ker}(L_n)=
\text{ker(ind}_{(G_{x_n},\Omega)}^{(G,\Omega)}(\tau_{{}_{f_n,x_n}},
\rho_{x_n}))$.

Pass to a sub-sequence, apply Lemma 8 to choose  
$$
\{f^\prime_n\}_{n=1}^\infty\subseteq\frak g^*, 
\{y_n\}_{n=1}^\infty\subseteq\Omega
$$
with
$f^\prime_n\to f$, $y_n\to x$, and
$$
\text{ker}(L_n)=\text{ker}
(\text{ind}^{(G,\Omega)}_{(G_{y_n},\Omega)}
(\tau_{{}_{f^\prime_n,y_n}}, \rho_{y_n})).
$$
 As $(f^\prime_n,y_n)\to (f,x)$ in $\go$, and $F$ is
saturated and closed, we are done.\qed

\vfill\eject

\topmatter
\title Chapter {3} \endtitle
\endtopmatter

\title Traces of irreducible representations of $C^*(G,\Omega)$\endtitle
\endtopmatter
\twoskip

We give a character theory for $C^*(G,\Omega)$ with
$G$ nilpotent Lie, analogous to Kirillov's character theory.
Our primary reference is Corwin and Greenleaf \cite{3}.

\twoskip

\noindent{\fod\bf Section 3.1}

\twoskip

\noindent{\fod\bf Schwartz functions on $G$}

\twoskip
\noindent{\bf Definition 17.}

In $\Bbb R^n$, define multi-indices $\alpha = (\alpha_1, \alpha_2, ... \alpha_n)$
and $\beta = (\beta_1, \beta_2, ... \beta_n)$. Assume
$\alpha_i\in\Bbb N_0^n$ and $\beta_i \in\Bbb N_0^n, \forall i$. 

Let $x^\beta=x_1^{\beta_1}\!\cdots x_n^{\beta_n},$
and $ \text{D}^\alpha=
\left(\frac{\partial}{\partial x_1}\right)^{\alpha_1}\!\cdots
\left(\frac{\partial}{\partial x_n}\right)^{\alpha_n}
$, a \lq\lq polynomial coefficient
differential operator" by $x^\beta D^\alpha$.

      On $\Bbb R^n$, the {\it Schwartz functions} $\Cal S(\Bbb R^n)$ are 
those
      $C^\infty$ functions $f$ such that
$$
\Vert x^\beta\text{D}^\alpha f\Vert_\infty<\infty,
$$
\noindent for all multi-indices $\alpha,\beta$. The 
natural topology of $\Cal S(\Bbb R^n)$ is determined by these
seminorms.
Denote the polynomial coefficient differential operators on 
$\Bbb R^n$
by $\Cal P(\Bbb R^n)$, and let $L\in\Cal P$ be arbitrary. 


    On $G$, define $\Cal S(G)$, the Schwartz functions on $G$, to 
be the $C^\infty$ functions on $G$ such that $\Vert L(f)\Vert_\infty$ is bounded for all
$L\in\Cal P(\Bbb R^n)$.

\twoskip

\noindent {\fod\bf Section 3.2}
\twoskip
\noindent{\fod\bf Traces on $C^*(G)$}
\twoskip

\smallskip

     For $\pi$ a unitary infinite-dimensional irreducible representation of a
nilpotent Lie group $G$, the operators $\pi(x)$ have no trace. However, for any
$\pi\in\widehat G$, there is a tempered distribution $\theta_\pi$ on $G$ that
plays the role of the classical trace character $\theta_\pi(g)=
\text{Tr}(\pi(g))$ for finite and compact groups. For $\pi\in\Cal S(G)$, the
Schwartz functions on $G$, the operator
$$
\pi(\phi)(\xi)=\int\limits_{s\in G}\phi(s)\pi(s)\,\xi\,\dd\mu_G(s)\ \ \
(\xi\in H_\pi)
$$
\noindent turns to be trace class; 
To wit:

\proclaim{Theorem 4}
Let $\pi=\pi_l$ be an irreducible representation of a
nilpotent Lie group, let $\frak m$ be a polarization for $l$, and model $\pi$ in
$L^2(\Bbb R^k)$ using any weak Malcev basis through $\frak m$. If
$\phi\in\Cal S(G)$, then $\pi(\phi)$ is trace class and
$$
\pi_\phi f(s)=\int_{\Bbb R^k} K_\phi(s,t)f(t)\,\dd\, t,\ \text{for all}\ 
f\in L^2(\Bbb R^k)
$$
\noindent where $K_\phi\in\Cal S(\Bbb R^k\times\Bbb R^k)$. Furthermore,
$$
\theta_\pi(\phi)=\text{Tr}(\pi(\phi))=\int_{\Bbb R^k}K_\phi(s,s)\dd s\ 
\text{(absolutely convergent)}
$$
\noindent and the functional $\theta_\pi$ is a tempered distribution on 
$\Cal S(G)$.
\endproclaim
\demo{Proof} See \cite{3}, Theorem 4.2.1, page 133.\qed

\twoskip
\noindent{\bf Definition 18.} 
     We may obtain explicit formulas for the kernel integral $K_\phi$,
once an $f\in\frak g^*$,
a polarization $\frak p$, and a weak Malcev basis
(\cite{3}, Theorem 1.1.13, pg. 10) 
through $\frak p$ are
specified. Let $n = \text{dim}(\frak g)$. Let $P=\text{exp}(\frak p)$; assume dim$(\frak g / \frak p) = k$. 
If $\{X_1,...,X_n\}$ is the weak Malcev basis, let
     $p=n-k=$ dim$(\frak p)$. Define polynomial maps
$\gamma:\Bbb R^n\mapsto G,\ \alpha:\Bbb R^p\mapsto P,\  \beta:
\Bbb R^k\mapsto G/P$ by
$$
\gamma(s,t)=\text{exp}(s_1X_1)\cdots\text{exp}(s_pX_p)\cdot
\text{exp}(t_1X_{p+1})\cdots\text{exp}(t_kX_n)
$$
$$
\alpha(s)=\gamma(s,0), \ \ \ \beta(t)=\gamma(0,t).
$$

 Let $\dd\mu_{{}_G},\ \dd\mu_{{}_P},\ \dd\mu_{{}_{G/P}}$
be the invariant measures on $G, P$ and $G/P$ determined by Lebesgue 
measures
$\dd s\,\dd t,\dd s,\ \dd t$, as in Theorems 1.2.10, 1.2.12
and 1.2.13 of \cite{3}.

\twoskip
    We describe Tr$(\pi(\phi))$ in terms
of integrals
over coadjoint orbits in $\frak g^*$. Given a 
Euclidean
measure $\dd X$ on $\frak g$, normalize measures on $\frak g$ and
$\frak g^*$ so that Fourier inversion holds, and define the
{\it Euclidean Fourier transforms}
$\widehat h$ (resp. $\Cal F h$) of functions $h$ on $G$
(resp. $\frak g$), as in \cite{3}, pg. 137.

 Each coadjoint orbit $\Cal O_{{}_f}=\text{Ad}^*(G)f$
has an invariant measure $\mu$ that is unique up to scalar multiple,
as
$\Cal O_{{}_f}\cong G/R_{{}_f}$, where
$$
R_{{}_f}=\text{Stab}_{{}_G}(f)=\{x\in G\ |\ \text{Ad}^*(x)f=f\}
$$

\twoskip

%
%

\noindent{\bf\fod\bf Section 3.3}
\twoskip
\noindent{\fod\bf Traces of irreducible representations of 
$C^*(G,\Omega)$}
\twoskip

Here we characterize some operators on
$C^*(G,\Omega)$
which are trace class. Assume $G$ is a connected, simply connected 
nilpotent Lie group with Lie algebra $\frak g$.

\smallskip

Recall that we use left actions.

\twoskip
\proclaim{Proposition 4} Let $\phi\in\Cal S(G)$, $\psi\in 
C_0(\Omega)$.
Let $L=(V,M)$ be an irreducible representation of $C^*(G,\Omega)$ corresponding to
$(f,x)\in\frak g^*\times\Omega$. Let $k=\text{dim}(\frak 
g/\frak p_x)$,
where $\frak p_x$ is polarizing for the action of $f$ in $\frak g_x$.
Assume the Hilbert space of $L$ is
$L^2(\Bbb R^k)$. For $\psi\in C_0(\Omega)$ with $\psi(\cdot\ x) |_{G/G_x}\in\Cal S(G/G_x)$
and $\phi\in\Cal 
S(G)$, then
$L(\phi\cdot\psi)$ has kernel $K$, 

$$
K(r,s)=\psi(\text{exp}(r)\cdot x)\int\limits_{t\in P_x}
\phi(\beta(r)t\beta(s)^{-1})e^{i\cdot f(\text{log}(t))}\dd\mu_{\frak p}(t)
$$. 

\endproclaim
\demo{Proof}\enddemo 

For $\phi\in\Cal S(G),\ \psi\in C_0(\Omega)$ such that 
$\psi(\cdot\ x) |_{G/G_x}\in\Cal S(G/G_x)$, $h\in L^2(\Bbb R^k)$ we have

$$
(L(\phi\cdot\psi)h)(\text{exp}(r))=\psi(\text{exp}(r)\cdot x)
\int\limits_{s\in \frak g}\phi(\text{exp}(s))h(\text{exp}(s)^{-1}\text{exp}(r))
\dd\mu_{\frak g}(s)=
$$

\twoskip
\centerline{(letting $s\to s^{-1}$, and our measure doesn't change)}
\twoskip

$$
\psi(\text{exp}(r)\cdot x)\int\limits_{s\in \frak g}
\phi(\text{exp}(s^{-1}))h(\text{exp}(s)\text{exp}(r))\dd\mu_{\frak g}(s)=
$$

\twoskip
\centerline{(letting $s\to sr^{-1}$)}
\twoskip

$$
\psi(\text{exp}(r)\cdot x)\int\limits_{s\in \frak g}
\phi(\text{exp}(r)(\text{exp}(s^{-1})))h(\text{exp}(s))\dd\mu_{\frak g}=
$$

\twoskip
\centerline{(splitting $\frak g$ into $\frak g/\frak p$ and $\frak p$:)}
\twoskip

$$
\psi(\text{exp}(r)\cdot x)
\int\limits_{s\in \frak g/\frak p}\int\limits_{t\in\frak p}
\phi(\text{exp}(r)(\text{exp}(t^{-1})\text{exp}(s^{-1})))
h(\text{exp}(s)\text{exp}(t))\dd\mu_{\frak g}=
$$
$$
\psi(\text{exp}(r)\cdot x)
\int\limits_{s\in \frak g/\frak p}\int\limits_{t\in\frak p}
\phi(\text{exp}(r)(\text{exp}(t^{-1})\text{exp}(s^{-1})))
e^{i\cdot f(t)}h(\text{exp}(s))\dd\mu_{\frak g}
$$

\noindent We now let
$$
K(r,s)=\psi(\text{exp}(r)\cdot x)\int\limits_{t\in\frak p}
\phi(\beta(r)\text{exp}(t)\beta(s)^{-1})e^{i\cdot f(t)}\dd\mu_{\frak p}(t),
$$
\noindent and the result is clear.\qed

\twoskip
\noindent{\bf Definition 19}

Let $H$ be a Hilbert space; for some $k\in\Bbb N$,
we assume that $H\cong L^2(\Bbb R^k)$.
For a
function $\psi\in C_0(\Bbb R^k)$, we denote the multiplication
operator by $\psi$ as
$M_\psi$.

\twoskip
\proclaim{Lemma 12} Let $\rho$ be an irreducible representation of a subgroup
$H$ of $G$, corresponding to the restriction of a functional $f\in\frak g^*$
to $\frak h$.
Defining $\pi=\text{ind}_H^G(\rho)$, we have that the the operator 
$T_{\psi,\phi}=M_\psi\cdot\pi(\phi)$ is trace class with a Schwartz kernel
when $\psi\in\Cal S(G/H)$ and $\phi\in\Cal S(G)$.
Also, for fixed $\psi$, $T_{\psi,\phi}$ is tempered in $\phi$. 
\endproclaim
\demo{Proof}\enddemo

     We do this by induction on the codimension of $H$ in $G$. If 
codim$(H)=0$, this is true by \cite{3}, Theorem 4.2.1. 

    Assume the lemma true for codim$(H)=n$.
   
By \cite{3}, Theorem 1.1.3 we find a subgroup
$G_1\subseteq G$ with codim$(G_1)=1$ and $H\subseteq G_1$.
The subgroup $H$ is codimension $n$ in $G_1$, and
$G_1$ is normal in $G$ by \cite{3}, Lemma 1.1.8.

 Assume the Lie algebra of $G_1$ is $\frak g_1$.
Let $\frak g=\frak g_1\oplus\bigl(\Bbb R-\text{span}\{X\}\bigr)$;
we have a smooth cross section for $G/G_1$ by exp$(\Bbb R\cdot X)\cong\Bbb R$.
  
    Define $\tau=\text{ind}_H^{G_1}(\rho)$, acting on the Hilbert space
$H_\tau$.

Assume $\pi=\text{ind}_H^G(\rho)\cong\text{ind}_{G_1}^G(\tau)$
acts on $H_\pi=L^2(\Bbb R)\otimes H_\tau$.

  Let $a\in G/G_1$, $f\in C^\infty(H_\pi)$.
By $f(a)$ we denote the element of $H_\tau$ corresponding to
$f(a)$. For $z\in G_1$, by $f(a)(z)$ we refer to the value of the
$H_\tau$-valued function $f(a)$ in $H_\tau$ at the point $z$ in $G_1$.

Let $\phi_1\in\Cal S(G/G_1)$,
$\phi_2\in\Cal S(G_1)$, $\psi_1\in\Cal S(G/G_1)$, and 
$\psi_2\in\Cal S(G_1/H)$. Note that sums of elements of the form
$\phi_1\cdot\phi_2$ are dense in $\Cal S(G)$, and that sums of the form
$\psi_1\cdot\psi_2$ are dense in $\Cal S(G/H)$.

 Let $r\in G/G_1$. Let $f\in C^\infty(H_\pi)$. We have
$$
(\rho_{\psi_1\cdot\psi_2}\pi(\phi_1\cdot\phi_2)f)(r)=
\rho_{\psi_1\cdot\psi_2}\int\limits_{s\in G/G_1}\int\limits_{t\in G_1}
\phi_1(s)\phi_2(t)\pi(s)\pi(t)f(r)\dd s\,\dd t=
$$
$$
\rho_{\psi_1\cdot\psi_2}
\int\limits_{s\in G/G_1}\int\limits_{t\in G_1}
\phi_1(s)\phi_2(t)\pi(s)f(t^{-1}r)\dd s\,\dd t=
$$
\centerline{(Note that $r^{-1}tr\in G_1$ by normality)}
\twoskip
$$
\rho_{\psi_1\cdot\psi_2}
\int\limits_{s\in G/G_1}\int\limits_{t\in G_1}
\phi_1(s)\phi_2(t)\pi(s)\tau(r^{-1}tr)f(r)\dd s\,\dd t=
$$
$$
\rho_{\psi_1\cdot\psi_2}\int\limits_{s\in G/G_1}\int\limits_{t\in G_1}
\phi_1(s)\phi_2(t)\tau(r^{-1}sts^{-1}r)f(s^{-1}r)\dd t\,\dd s=
$$
\twoskip
\centerline{(Let $t\to s^{-1}rtr^{-1}s$ in the inner integral and set $s^{-1}r = a$)}
\twoskip
$$
\rho_{\psi_1}\int\limits_{a\in G/G_1}\phi_1(ra^{-1})\left[\rho_{\psi_2}
\int\limits_{t\in G_1}\phi_2(ata^{-1})\tau(t)\dd t\right] f(a)\dd a.\tag1 
$$

\noindent Note that for fixed $a$, by induction hypothesis
the operator in the brackets of (1) above is a trace class operator on $H_\tau$,
and for fixed $\psi_2$, is tempered in $\phi$. 
Employ \cite{3}, Proposition 1.2.8 and for some selection of polynomials
$\{P_i\}_{i=1}^n$ and at a fixed point $y\in G_1$ the inner integral of 
formula (1) equals
$$
\multline
\left(\rho_{\psi_2}
\int\limits_{t\in G_1}
\phi_2\left( P_1(a,t),\dots P_n(a,t)\right)\tau(t)f(a)\dd t\right)(y)= \\
\int\limits_{z\in G_1}k^{a}_{\phi_2}(z,y)f(a)(z)\dd z.
\endmultline
$$
\noindent By the inductive hypothesis,
  we may find the integral kernel $k^{a}_{\phi_2}(z,y)$ 
 and it may be chosen Schwartz in $y$ and $z$.

Recall that for fixed $a$, Tr$(k_{\phi_2}^{a})$ is 
tempered in $\phi_2$ by inductive hypothesis. By \cite{27}, Corollary 1, pg. 43
we know that  there 
exists a seminorm $p$ on $G_1$ such and a constant $C$ such that
that $|\text{Tr}(k_{\phi_2}^{a}))|\le C\cdot p(\phi_2^{a})$,
where $\phi_2^{a}$ is defined as obvious.

    Assume that we have  multi-indices $\alpha,\beta\in\Bbb N^n_0,\ \text{with}\ 
\alpha=(\alpha_1,\dots,\alpha_n),$\linebreak $\beta=\beta_1,\dots,\beta_n$;
define $L = 
 w_1^{\beta_1}\cdots w_n^{\beta_n}\cdot
\frac{\partial^{|\alpha |}}{\partial w_1^{\alpha_1}\cdots \partial w_n^{\alpha_{i-1}}}$

And we have for the seminorm $p$:
$$
|p(\phi^a)|= 
\underset{w\in\Bbb R^n}\to{\text{sup}}
\left\{\left |\
L\bigl(\phi_2\left(P_1(a,w),\dots,P_n(a,w)\right)\bigr)\right |\ 
\right\}\le 
|Q(a)|,
$$

\noindent by Proposition 1.2.9 of \cite{3}
and the $n$ dimensional chain rule, where $Q$ is some polynomial.

So Tr$(k_{\phi_2}^a)=\int_{y\in G_1}k_{\phi_2}^a(y,y)\dd y$ 
grows no faster than polynomial in $a$, and $k^a_{\phi_2}$ grows no 
faster
that polynomial in $a$.

For $a$ in a bounded set, $\phi_2(ata^{-1})$ is
bounded by a $L^1$ function. 
   Let $\frak p_x$ be a polarizing subalgebra of $\frak g_x$ with
   respect to the restriction of $f$ to $\frak g_x$.
 By Proposition 4 we have:
$$
k^a_{\phi_2}(r,s)=\psi_2(\text{exp}(r)\cdot x)\int\limits_{t\in P_x}
\phi_2(\beta(r)\cdot ata^{-1}\cdot\beta(s)^{-1})e^{i\cdot f(\text{log}(t))}
\dd\mu_{\frak p}(t),
$$
\noindent and we may differentiate
in $a$ under the integral sign; infinite
differentiability of $k_{\phi_2}^a$ in $a$ follows.

   Choose a smooth splitting of $G$ into
$G/G_1$ and $G_1$ by projections $p:G\mapsto G/G_1$ and $q:G\mapsto G_1$.

    For our functions $\phi_1,\phi_2,\psi_1$ and $\psi_2$,
 define the integral kernel $K$ of \linebreak
$
\rho_{\psi_1\cdot\psi_2}\cdot\pi(\phi_1\cdot\phi_2)$ by
$$
K(x,y)=\rho_{\psi_1}\bigl(p(y)\bigr)\cdot\phi_1\bigl(p(y)p(x^{-1})\bigr)
\cdot k^{p(y)}_{\phi_2}\bigl(q(x),q(y)\bigr).
$$
\noindent Integrated against $f\in H_\pi$, this gives us
$\rho_{\psi_1\cdot\psi_2}\cdot\pi(\phi_1\cdot\phi_2)(f)$:

$$
(\rho_{\psi_1\cdot\psi_2}\cdot\pi(\phi_1\cdot\phi_2)f)(y) =
(Kf)(y)=\int\limits_{x\in G}K(x,y)f(x)\dd x=
$$
$$ 
\rho_{\psi_1}(p(y)) 
\int\limits_{x\in G}\phi_1(p(y)p(x^{-1}))k^{p(y)}_{\phi_2}
(q(x),q(y))f(x)\dd x=
$$
$$
\rho_{\psi_1}(p(y))\int\limits_{a\in G/G_1}\int\limits_{z\in G_1}
\phi_1(p(y)a^{-1})k^a_{\phi_2}(z, q(y))f(a)(z)\dd z\,\dd a.
$$

\noindent Note $k_{\phi_2}^a$ is Schwartz on $G_1$, and its
integral and the integrals of its derivatives grow no faster than
polynomial in $a$.

   Note $G/G_1\cong\Bbb R$; treat 
composition on
$G/G_1$ as addition on $\Bbb R$.
The functions $\psi_1$ and $\phi_1$ are both Schwartz on
$G/G_1$, and part of the integral kernel is
$\rho_{\psi_1}\bigl(p(y)\bigr)\cdot\phi_1\bigl(p(y)p(x^{-1})\bigr)$.
The other part $k^a_{\phi_2}$ is Schwartz already.

By Peetre's Inquality (\cite{14}, pg. 10), when $\psi$ and $\phi$ are
both Schwartz on $\Bbb R$,
$\psi(y)\phi(y+a)$ is Schwartz on $\Bbb R\times\Bbb R$.

Similar properties follow for derivitives,
and $\psi(y)\phi(y+a)$ is Schwartz on $\Bbb R^2$.

    Our entire integral kernel is Schwartz, and the operator 
$\rho_{\psi_1\cdot\psi_2}\cdot\pi(\phi_1\cdot\phi_2)$ is trace class by
\cite{3} Theorem A.3.9.

The kernel $K$ is clearly tempered in $\phi_1\in\Cal S(G/G_1)$;
the final result follows by induction and density arguments. \qed
\twoskip

Henceforth assume that $\Omega/G$ is a $T_0$ space.
By Theorem 2.1 \cite{8}, $(G,\Omega)$ is Polish, and
for any $x\in\Omega$, we have $G\cdot x\cong G/G_x$. 

\smallskip
\noindent{\bf Definition 20.}
     For any orbit $G\cdot x$, define
$$
\Cal A_{{}_{G\cdot x}}=\Bbb R-\text{span}\bigl\{\phi\cdot\psi\ |\
\phi\in\Cal S(G),\
\psi\in C_0(\Omega)\ \text{and}\ \psi(\cdot\ x)|_{G/G_x}\in
 \Cal S(G/G_x)\bigr\}.
$$

\twoskip
\proclaim{Theorem 5} Let $L=(V,M)$ be the irreducible
representation of $C^*(G,\Omega)$
associated to the pair
$(f,x)\in\frak g^*\times\Omega$. The representation $L$ is trace class on
$\Cal A_{{}_{G\cdot x}}$, and for fixed $\psi$, is tempered in 
$\phi$.
\endproclaim
\demo{Proof}\enddemo Note
$L(\psi\cdot\phi)=M_\psi(\cdot\ x)\cdot V(\phi)$, and
apply Lemma 12 just proven.\qed

\twoskip
\noindent{\bf Definition 21.}
Remember Definition 18, where we defined maps $\alpha,\ \beta,\ 
\gamma$ of $\frak p$, $\frak g/\frak p$, and $\frak g$ (resp.)
to $G$. Assume $f\in\frak g^*$, and for $x\in\Omega$,
$\frak g_x$ is the Lie algebra of $G_x$.
Assume dim$(\frak g_x)=l$. Define a new map
$\delta :\Bbb R^l\mapsto G_x$ by
$\delta(s)=\gamma(s_1,\dots,s_l,0)$.

\twoskip
\proclaim{Theorem 6} Let $L$ be a representation 
of $C^*(G,\Omega)$
corresponding to the functional-point pair $(f,x)\in\go$.

Let $\frak p_x$ be a polarizing subalgebra of $\frak g_x$ 
with respect to $f|_{\frak g_x}$;
assume $k=\text{dim}(\frak g/\frak p_x)$.

Assume $L$ acts on the Hilbert space
$L^2(\Bbb R^k)$.

Let $\Cal O_{{}_L}$ be the 
$\sim_1$ equivalence class of $(f,x)$ in $\go$, specifically,
$$
\Cal O_{{}_{(f,x)}}=\left\{ (l,y)\in\go\ 
|\ \text{for some}\ s\in G,
\text{we have} \right.
$$
$$
\left. l=\text{Ad}^*(s)f+h,\
y=s\cdot x,
\ h\in\frak g_x^\perp\right\}.
$$

Let $p,q$ be the natural projections from $\Cal O_{{}_L}$
to $\frak g^*$ and $\Omega$, respectively. Let 
$\phi\in\Cal S(G)$,
$\psi(\cdot\ x)\in\Cal S(G/G_x)$. We have
$$
\text{Tr}(L(\phi\cdot\psi)
)=\int\limits_{z\in\Cal O_{{}_L}}
\psi(q(z))\widehat\phi(p(z))\dd z
$$
\noindent for a particular choice of $G$-invariant 
measure $\dd z$
on $\Cal O_{{}_L}$.
\endproclaim
\demo{Proof}\enddemo

    We closely mimic the proof of Theorem 4.2.4 on pp. 138-41 of \cite{3}.

    Assume that $\frak p_x$ has dimension $p$, so $n=\text{dim}(\frak g)=p+k$.

     We give $\frak g$ a
standard basis realization on $L^2(\Bbb R^k)$. By Proposition 4 
above, and \cite{3}, Theorem A.3.9, we have
$$
\text{Tr}(L)=\int\limits_{s\in\frak g/\frak p_x}K(s,s)\dd s=
$$
$$
\int\limits_{s\in\frak g/\frak p_x}\psi(\beta(s)\cdot x)
\int\limits_{t\in\frak p_x}\phi(\beta(s)\text{exp}(t)\beta(s)^{-1})e^{i\cdot f(t)}
\dd\mu_{\frak p_x}(t)\,\dd\mu_{\frak g/\frak p_x}(s).\tag1
$$

    Let $\frak z$ be a subspace complementary to $\frak p_x$ in $\frak g$. It is
easy to take $\frak z=\Bbb R$-span of the last $k$ basis vectors of our
Malcev basis through $\frak p_x$. We have an additive splitting $H+X\in\frak p_x
\oplus\frak z$ for each element in $\frak g$. Let $\dd X,\ \dd H$ be arbitrarily
assigned Euclidean measures on $\frak z$, $\frak p_x$; then we have that
$\dd H\,\dd X$ is a Euclidean measure on $\frak g$, which we use to define
the above integrals. 

     For $\phi\in\Cal S(G),\ u\in\Bbb R^k\cong\frak z\cong\frak g_x/\frak p_x$, 
we define
$$
\phi_u(H,X)=\phi(\beta(u)\text{exp}(H+X)\beta(u)^{-1}).
$$
\noindent For each fixed $u$, this is a Schwartz function on $\Bbb R^p\times
\Bbb R^k\cong\frak p_x\times\frak z$. Viewing $\frak z^*\cong\frak p_x^\perp$,
define the Fourier transform
$$
\Cal F \phi(f^\perp)=\int\limits_{X\in\frak z}\phi(X)e^{i\cdot f^\perp(X)}
\dd f^\perp\ \ \ \ \text{for}\ \phi\in\Cal S(Z), \ f^\perp\in\frak p_x^\perp.
$$
\noindent When the measures are suitably normalized, this is an $L^2$
isometry with inverse
$$
(\Cal F^{-1}\phi)(X)=\int\limits_{f^\perp\in\frak p_x^\perp}\phi(f^\perp)
e^{-i\cdot f^\perp(X)}\dd f^\perp,
$$
\noindent taking the dual Euclidean measure $\dd l^\perp$ on $\frak p_x^\perp$.
We have
$$
\phi_u(H,0)=(\Cal F^{-1}\Cal F\phi_u)(H,0)=
\int\limits_{f^\perp\in\frak z^*}
e^{-i\cdot f^\perp(0)}\Cal F\phi_u(H,f^\perp)\dd f^\perp=
$$
$$
\int\limits_{f^\perp\in\frak p_x^\perp}\int\limits_{X\in\frak z}
\phi_u(H,X)e^{i\cdot f^\perp(X)}\dd X\,\dd f^\perp \tag2
$$
 
\noindent We want to insert equation 1 into equation
2 and interchange some integrals, this 
will
be done introducing and removing an ad hoc function which enables us to use
Fubini's Theorem. Let $\{w_j\}_{j=1}^\infty$ be a collection of functions on
$\frak p_x^\perp$ with $0\le w_j(f^\perp)\le 1$, and $w_j\uparrow 1$ uniformly on
compacta in $\frak p_x^\perp$ and for all $j$,

$$
\int\limits_{f^\perp\in\frak p_x^\perp}w_j(f^\perp)\dd f^\perp<\infty
$$

As $\Cal F\phi_u(H,f^\perp)$ is Schwartz in both variables (remember that the 
Fourier transform of a Schwartz function is also Schwartz; see Proposition 
4 of \cite{27}), by dominated convergence we have
$$
\text{Tr}(L)=\int\limits_{u\in\frak g/\frak p_x}\psi(\beta(u)\cdot x)
\int\limits_{H\in\frak p_x}e^{i\cdot f(H)}\phi_u(H,0)\dd H\,\dd u=
$$
$$
\int\limits_{u\in\frak g/\frak p_x}\psi(\beta(u)\cdot x)\int\limits_{H\in\frak p_x}
\biggl[\lim_{j\to\infty}\int\limits_{f^\perp\in\frak p_x^\perp}
\Cal F\phi_u(H,f^\perp)w_j(f^\perp)\dd f^\perp\biggr]e^{i\cdot f(H)}
\dd H\,\dd u=
$$
$$
\multline
\int\limits_{u\in\frak g/\frak p_x}\psi(\beta(u)\cdot x)\biggl[\lim_{j\to\infty}
\int\limits_{H\in\frak p_x}\int\limits_{f^\perp\in\frak p_x^\perp}
\int\limits_{X\in\frak z}e^{i\cdot f(H)}e^{i\cdot f^\perp(X)}
\phi_u(H,X)w_j(f^\perp)\cdot \\
\dd X\,\dd f^\perp\,\dd H\biggr]\dd u.
\endmultline
$$

 Fubini's Theorem may be applied to the innermost triple integral, when
we re-arrange, we get
$$
\text{Tr}(L)=
$$
$$
\multline
\int\limits_{u\in\frak g/\frak p_x}\psi(\beta(u)\cdot x)\biggl[
\lim_{j\to\infty}\int\limits_{f^\perp\in\frak p_x^\perp}w_j(f^\perp)
\int\limits_{\{H,X\}\in\frak p_x\oplus\frak z}
e^{i\cdot\left(f(H)+f^\perp(X)\right)}\phi_u(H,X) \\
\dd X\,\dd H\,\dd f^\perp\biggr]\dd u=
\endmultline
$$
$$
\int\limits_{u\in\frak g/\frak p_x}\psi(\beta(u)\cdot x)\biggl[
\int\limits_{f^\perp\in\frak p_x^\perp}\
\int\limits_{\{H,X\}\in\frak p_x\oplus\frak z}
e^{i\cdot\left(f(H)+f^\perp(X)\right)}\phi_u(H,X)
\dd X\,\dd H\,\dd f^\perp\biggr]\dd u,
$$
\noindent by the Dominated Convergence Theorem, as the integral over 
$\frak p_x\oplus\frak z$ is Schwartz in $f^\perp$. The integral over 
$\frak p_x^\perp$ amounts to integration over $\frak z^*$; by translation
invariance in this integral, we may replace $f^\perp(X)$ with $(f+f^\perp)(X)$.
As $f(H)=(f+f^\perp)(H)$ for all $f^\perp\in\frak p_x^\perp$, all $H\in\frak p_x$,
we may write
$$
\text{Tr}(L)=
$$
$$
\multline
\int\limits_{u\in\frak g/\frak p_x}\psi(\beta(u)\cdot x)
\int\limits_{f^\perp\in\frak p_x^\perp}
\int\limits_{\{H,X\}\in\frak p_x\oplus\frak z}
e^{i\cdot(f+f^\perp)(H)}e^{i\cdot(f+f^\perp)(X)} \\
\phi\bigl(\beta(u)\text{exp}(H+X)\beta(u)^{-1}\bigr)
\dd X\,\dd H\,\dd f^\perp\,\dd u=
\endmultline
$$
\twoskip
$$
\int\limits_{u\in\frak g/\frak p_x}\psi(\beta(u)\cdot x)
\int\limits_{f^\perp\in\frak p_x^\perp}\int\limits_{Y\in\frak g}
e^{i\cdot(f+f^\perp)(Y)}
\phi\bigl(\beta(u)\text{exp}(Y))\beta(u)^{-1}\bigr)
\dd Y\,\dd f^\perp\,\dd u=
$$
\twoskip
\centerline{(letting $Y\to\text{Ad}(\beta(u^{-1}))(Y)$ 
in the last integral above)}
\twoskip
$$
\int\limits_{u\in\frak g/\frak p_x}\psi(\beta(u)\cdot x)
\int\limits_{f^\perp\in\frak p_x^\perp}\int\limits_{Y\in\frak g}
e^{i\cdot(f+f^\perp)(\text{Ad}(\beta(u^{-1}))(Y))}
\phi(\text{exp}(Y)) 
\dd Y\,\dd f^\perp\,\dd u=
$$
$$
\int\limits_{u\in\frak g/\frak p_x}\psi(\beta(u)\cdot x)
\int\limits_{f^\perp\in\frak p_x^\perp}\int\limits_{Y\in\frak g}
e^{i\cdot(\text{Ad}^*(\beta(u))(Y)f+f^\perp)}
\phi(\text{exp}(Y)) 
\dd Y\,\dd f^\perp\,\dd u=
$$
$$
\int\limits_{u\in\frak g/\frak p_x}\psi(\beta(u)\cdot x)
\int\limits_{f^\perp\in\frak p_x^\perp}
\widehat\phi\bigl(\text{Ad}^*(\beta(u)^{-1})(f+f^\perp)\bigr)
\dd f^\perp\,\dd u.\tag3
$$
\twoskip

   We now split $\frak g/\frak p_x$ into $\frak g/\frak g_x$ and $\frak g_x/\frak p_x$,
and as $\beta(t)\cdot x=x$ for all $t\in\frak g_x/\frak p_x$, and we get that
formula (3) above equals
$$
\int\limits_{u\in\frak g/\frak g_x}\psi(\beta(u)\cdot x)
\int\limits_{v\in\frak g_x/\frak p_x}\int\limits_{f^\perp\in\frak p_x^\perp}
\widehat\phi\bigl(\text{Ad}^*(\beta(u))\text{Ad}^*(\beta(v))(f+f^\perp)\bigr)
\dd f^\perp\,\dd v\,\dd u=
$$
\twoskip
$$
\multline
\int\limits_{u\in\frak g/\frak g_x}\psi(\beta(u)\cdot x)
\int\limits_{v\in\frak g_x/\frak p_x}
\int\limits_{f_2^\perp\in\frak g_x^\perp}
\int\limits_{f_1^\perp\in\frak p_x^\perp/\frak g_x^\perp}
\widehat\phi\bigl(\text{Ad}^*(\beta(u))\text{Ad}^*(\beta(v))
(f+f_1^\perp+f_2^\perp)\bigr)
\\
\dd f_1^\perp\,\dd f_2^\perp\,\dd v\,\dd u.
\endmultline\tag4
$$
\twoskip

   We note that by using the integral in $f_2^\perp$ that we may now assume that
$f$ is nonzero only on $\frak g_x$.

     We now work with the integral in $f_1^\perp$. If $R_f$ is the
stabilizer of the functional $f$ restricted to $\frak g_x$, there exists 
invariant measures
$\dd\dot p$ and $\dd\dot x$ on $P_x/R_f$ and $G_x/P_x$ such that
$\dd\dot p\, \dd\dot x$
is invariant measure on $G_x/R_f$. We know from Proposition 3.1.18 of \cite{3} that 
$\text{Ad}^*(P_x)(f)=f+\frak p_x^\perp/\frak g_x^\perp
=(f+\frak p_x^\perp)|_{\frak g_x}$,
and that the natural diffeomorphism
$\Delta : P_x/R_f\mapsto\text{Ad}^*(P_x)(f)=
(f+\frak p_x^\perp)|_{\frak g_x}$ is
equivariant and measure preserving on $(f+\frak p_x^\perp)_{\frak g_x}$. Using
different ideas, we may also transfer Euclidean measure
$\dd f_1^\perp$ on $\frak p_x^\perp/\frak g_x^\perp$ under the translation map $q$,
where $q(f_1^\perp)=f+f_1^\perp$, to a Euclidean measure $\nu=q^*(\dd f_1^\perp)$
on the affine space $f+\frak p_x^\perp/\frak g_x^\perp$. Now for each $p\in P_x$, 
we define an affine map $A(p)$ from $\frak p_x^\perp/\frak g_x^\perp$ to itself by
$A(p)(f_1^\perp)=\text{Ad}^*(p)f_1^\perp+(\text{Ad}^*(p)f-f)$,
we note that
$(\text{Ad}^*(p)f-f)\in\frak p_x^\perp/\frak g_x^\perp$; see
Proposition 3.1.18 of \cite{3}. We now note that
$$
q\circ A(p)(f_1^\perp)=
q\bigl(\text{Ad}^*(p)(f_1^\perp)+(\text{Ad}^*(p)f-f)\bigr)=
\text{Ad}^*(p)(f+f_1^\perp).\tag5
$$
\noindent As the linear part $\text{Ad}^*(p)|_{\frak p_x^\perp/\frak g_x^\perp}$
of $A(p)$ is unipotent, the operator $A(p)$ preserves $\dd f_1^\perp$,
and by formula (1) above, $\text{Ad}^*(p)$ preserves $\nu$ on the affine space 
$f+\frak p_x^\perp/\frak g_x^\perp$. As $\text{Ad}^*(p)$ is also measure preserving 
on $P_x/R_f$, we have that under $\Delta$, $\nu$ is identified with an
invariant measure on $P_x/R_f$, which must be a scalar multiple of $\dd\dot p$:
$(\Delta^{-1})^*\nu=c\cdot\dd\dot p$. Hence if
$\phi\in L^1(f+\frak p_x^\perp/\frak g_x^\perp)$ we must have
$$
\int\limits_{f_1^\perp\in\frak p_x^\perp/\frak g_x^\perp}\phi(f+f_1^\perp)
\dd f_1^\perp=
\int\limits_{f^\prime\in f+\frak p_x^\perp/\frak g_x^\perp}\phi(f^\prime)
\dd\nu(f^\prime)=
c\cdot\int\limits_{\dot p\in P_x/R_f}\phi(\text{Ad}^*(\dot p)f)\dd\dot p.
$$ 
\twoskip

     Consequently, formula (4) above equals

\twoskip
$$
\multline
\int\limits_{u\in\frak g/\frak g_x}\psi(\beta(u)\cdot x)
\int\limits_{v\in\frak g_x/\frak p_x}
\int\limits_{f_2^\perp\in\frak g_x^\perp}
c\cdot\int\limits_{m\in\frak p_x/\frak r_f} \\
\widehat\phi\bigl(\text{Ad}^*(\beta(u))\text{Ad}^*(\delta(v))
(\text{Ad}^*(\alpha(m)) f+f_2^\perp)\bigr) \\
\dd m\,\dd f_2^\perp\,\dd v\,\dd u=
\endmultline
$$
\twoskip
\centerline{(by Haar measure we move the $\text{Ad}^*(\alpha(m))$ to include}
\centerline{$f_2^\perp$, we note that $\text{Ad}^*(\alpha(m))f_2^\perp$
is still zero on $\frak g_x$)}
\twoskip
$$
\multline
c\cdot\int\limits_{u\in\frak g/\frak g_x}\psi(\beta(u)\cdot x) \\
\int\limits_{v\in\frak g_x/\frak p_x}
\int\limits_{f_2^\perp\in\frak g_x^\perp}
\int\limits_{m\in\frak p_x/\frak r_f}
\widehat\phi\bigl(\text{Ad}^*(\beta(u))\text{Ad}^*(\delta(v))
(\text{Ad}^*(\alpha(m))(f+f_2^\perp))\bigr)
 \\
\dd m\,\dd f_2^\perp\,\dd v\,\dd u=
\endmultline
$$
\twoskip
\centerline{(now combining $v$ and $m$ into a single variable $y$)}
\twoskip

$$
c\cdot\int\limits_{u\in\frak g/\frak g_x}\psi(\beta(u)\cdot x)
\int\limits_{f_2^\perp\in\frak g_x^\perp}
\int\limits_{y\in\frak g_x/\frak r_{{}_f}}
\widehat\phi\bigl(\text{Ad}^*(\beta(u))\text{Ad}^*(\delta(y))(f+f_2^\perp)\bigr)
\dd y\,\dd f_2\,\dd u.
$$
\twoskip

     This is our orbital integral.\qed

\vfill\eject

\end